# ASYMPTOTIC OPTIMALITY OF MAXIMUM PRESSURE POLICIES IN STOCHASTIC PROCESSING NETWORKS[1]

By J. G. Dai and Wuqin Lin

*Georgia Institute of Technology and Northwestern University*

We consider a class of stochastic processing networks. Assume that the networks satisfy a complete resource pooling condition. We prove that *each* maximum pressure policy asymptotically minimizes the workload process in a stochastic processing network in heavy traffic. We also show that, under each quadratic holding cost structure, there is a maximum pressure policy that asymptotically minimizes the holding cost. A key to the optimality proofs is to prove a state space collapse result and a heavy traffic limit theorem for the network processes under a maximum pressure policy. We extend a framework of Bramson [*Queueing Systems Theory Appl.* **30** (1998) 89–148] and Williams [*Queueing Systems Theory Appl.* **30** (1998b) 5–25] from the multiclass queueing network setting to the stochastic processing network setting to prove the state space collapse result and the heavy traffic limit theorem. The extension can be adapted to other studies of stochastic processing networks.

**1. Introduction.** This paper is a continuation of Dai and Lin (2005), in which maximum pressure policies are shown to be throughput optimal for a class of stochastic processing networks. Throughput optimality is an important, first-order objective for many networks, but it ignores some key secondary performance measures like queueing delays experienced by jobs in these networks. In this paper we show that maximum pressure policies enjoy additional optimality properties; they are asymptotically optimal in

Received July 2006; revised February 2008.

[1]Research supported in part by NSF Grants CMMI-0727400 and CNS-0718701, and by an IBM Faculty Award.

*AMS 2000 subject classifications.* Primary 90B15, 60K25; secondary 90B18, 90B22, 68M10, 60J60.

*Key words and phrases.* Stochastic processing networks, maximum pressure policies, backpressure policies, heavy traffic, Brownian models, diffusion limits, state space collapse, asymptotic optimality.

minimizing a certain workload or holding cost of a stochastic processing network.

Stochastic processing networks have been introduced in a series of three papers by Harrison (2000, 2002, 2003). In Dai and Lin (2005) and this paper we consider a special class of Harrison's model. This class of stochastic processing networks is much more general than multiclass queueing networks that have been a subject of intensive study in the last 20 years; see, for example, Harrison (1988), Williams (1996) and Chen and Yao (2001). The added features in a stochastic processing network allow one to model skills-based routing in call centers [Gans, Koole and Mandelbaum (2003)], operator-machine interactions in semiconductor wafer fabrication facilities [Kumar (1993)], and combined input- and output-queued data switches [Chuang et al. (1999)] in data networks.

For this general class of stochastic processing networks, Dai and Lin (2005) propose a family of operational policies called maximum pressure policies and prove that they are throughput optimal. For a given vector $\alpha > 0$, the maximum pressure policy associated with the parameter $\alpha$ is specified in Section 3 of Dai and Lin (2005) and will be specified again in Definition 1 of this paper. In this paper, for a stochastic processing network that satisfies a complete resource pooling condition and a heavy traffic condition, we first show in Theorem 2 that a certain workload process is asymptotically minimized under *any* maximum pressure policy. When the holding cost rate is a quadratic function of the buffer contents, we show then in Theorem 3 that there is a maximum pressure policy that asymptotically minimizes the holding cost. Our maximum pressure policies do not solve the linear holding cost optimization problem for stochastic processing networks. However, following an approach in Stolyar (2004), one can find a maximum pressure policy that is *asymptotically* $\varepsilon$*-optimal* under the linear holding cost structure. Section 9 elaborates the $\varepsilon$-optimality of maximum pressure policies.

In Theorem 2, except for some nonnegativity requirements, the parameter $\alpha$ that is used to define a maximum pressure policy can be chosen arbitrarily. In Theorem 3, to minimize a given quadratic holding cost rate, one has to choose the parameter $\alpha$ to be the vector of coefficients that define the quadratic holding cost rate function. In both cases, the parameter can be chosen to be independent of network data like arrival or processing rates. Thus, these maximum pressure policies do not depend on the arrival rates. This feature is attractive in some applications when network data like arrival rates are sometimes difficult or impossible to be estimated accurately. The maximum pressure policies in Section 9 do depend on the arrival rates of the network.

Our asymptotic region is when the stochastic processing network is in *heavy traffic*; at least one server has to be 100% busy in order to handle all the input. A key assumption on our network is that a *complete resource pooling* condition is satisfied. Roughly speaking, the complete resource pooling condition requires enough overlap in the processing capabilities of bottleneck servers that these servers form a single, pooled resource or "super server." As will be discussed fully in Section 3, the complete resource pooling condition is articulated by the dual problem of a linear program (LP) called the *static planning problem*. For a network satisfying the complete resource pooling condition, the corresponding dual LP has a unique optimal solution and the workload process is defined by this unique optimal solution. Ata and Kumar (2005) develop a discrete-review policy and prove its asymptotic optimality in minimizing the linear holding cost for a class of *unitary* stochastic processing networks that satisfy the complete resource pooling condition and the *balanced* heavy traffic condition. The latter condition requires every server in the network to be heavily loaded. This balanced load requirement, combined with the complete resource pooling assumption, rules out some well-known networks such as multiclass queueing networks, which have been used to model semiconductor fabrication lines [Kumar (1993)]. Our definition of heavy traffic condition is less restrictive than those in Ata and Kumar (2005); our stochastic processing networks include those multiclass queueing networks that have a unique bottleneck server. Note that the nonbottleneck stations may not disappear in a heavy traffic diffusion limit in a multiclass queueing network operating under a nonidling service policy such as first-come–first-serve; see, for example, Bramson (1994). Even under an asymptotically optimal maximum pressure policy studied in this paper, the non-bottleneck stations may not disappear; see the example below Theorem 4 in Section 5.

The major part of our optimality proof of the maximum pressure policies is a heavy traffic limit theorem. The theorem asserts that when the network is operated under a maximum pressure policy, (a) the one-dimensional workload process converges to a reflecting Brownian motion in diffusion limit, and (b) the multidimensional buffer content process is a constant multiple of the workload process in diffusion limit. The latter result is a form of *state space collapse* for network processes, and its proof is the key to the proof of the limit theorem. We choose to extend a framework of Bramson (1998) and Williams (1998b), from the multiclass queueing network setting to the stochastic processing network setting, to prove the heavy traffic limit theorem. We first show that all solutions to a critically loaded fluid model operating under a maximum pressure policy exhibit some form of state space collapse. Then we translate the state space collapse to the diffusion scaling following Bramson (1998), proving the state space collapse in diffusion limit. Once we have the state space collapse result, we invoke a theorem of Williams (1998b) for perturbed Skorohod problems to establish the heavy traffic limit theorem for our stochastic processing network operating under the maximum pressure policy.

Stolyar (2004) proves that MaxWeight policies asymptotically minimize the workload processes in heavy traffic for a generalized switch model that belongs to one-pass systems in which each job leaves the system after being processed at one processing step. Our Theorem 2 greatly generalizes Stolyar (2004) from one-pass systems to stochastic processing networks. Except for Ata and Kumar (2005) that was discussed earlier in this introduction, most other works that are closely related to our work have focused on parallel server systems. These systems belong to a special class of one-pass systems. All these works assume heavy traffic and complete resource pooling conditions. For a parallel server system that has 2 buffers, 2 processors and 3 activities, Harrison (1998) develops a "discrete-review" policy via the BIG-STEP procedure that was first described in Harrison (1996) for multiclass queueing networks. Furthermore, he proves that the discrete-review policy asymptotically minimizes the expected discounted, cumulative linear holding cost under the restrictive assumption of Poisson arrival processes and deterministic service times. Harrison and López (1999) then use the BIG-STEP procedure to produce a family of policies for general parallel systems, but they have not proved the optimality of these policies. For the same 2 buffer, 3 processor and 3 activity parallel server system, but with general arrival processes and service time distributions, Bell and Williams (2001) develop simple form, buffer priority policies with thresholds and prove their asymptotic optimality under linear holding cost. Since the threshold values are constantly monitored, these policies are termed as "continuous-review" policies. They further generalize their policies to general parallel server systems in Bell and Williams (2005) and prove that they are asymptotically optimal. While all these works deal with the holding cost objective, the proposed asymptotically optimal policies in the literature exploit the special network structures and critically depend on the network data, particularly the arrival rates. Mandelbaum and Stolyar (2004) propose a generalized $c\mu$ policy for parallel server systems. The policy does not use any arrival rate information. They prove that it is asymptotically optimal in minimizing a strictly convex holding cost.

When a network has multiple bottlenecks so that the complete resource pooling condition is not satisfied, finding an asymptotically optimal policy remains a difficult, open problem in general. Shah and Wischik (2006) study the fluid and diffusion limits under MaxWeight policies for input-queued switches that do not satisfy the complete resource pooling condition. They propose a policy that is believed to be asymptotically optimal. However, they do not provide a proof for the optimality. Harrison and Wein (1989) study a two-station multiclass queueing network known as crisscross network. They propose a threshold type policy and demonstrate its near-optimal performance through simulations. For the same crisscross network, but with exponentially distributed interarrival and service times, Martins, Shreve and Soner (1996) and Budhiraja and Ghosh (2005) prove the asymptotic optimality of certain nested threshold policies when the network data is in various heavy traffic regimes.

All asymptotic optimality proofs in the literature involve proving a heavy traffic limit theorem and some form of state space collapse, either explicitly or implicitly. Ata and Kumar (2005) and Bell and Williams (2001, 2005) prove the state space collapse results directly without going through fluid models. Stolyar (2004) and Mandelbaum and Stolyar (2004) mimic the general framework of Bramson (1998) and Williams (1998b). They start with showing a state space collapse result for fluid models, and then prove the optimality directly without proving the state space collapse in diffusion limit as an intermediate step. By choosing to extend Bramson and Williams' framework in this paper, we are able to provide an optimality proof that is clean and hopefully easy to follow. We expect our extension can be adapted to other studies of stochastic processing networks. Our proof of asymptotic optimality requires that the service times have finite $2+\varepsilon$ moments, as in Ata and Kumar (2005). This moment assumption is weaker than the exponential moment assumption that is usually assumed in the literature; see, for example, Harrison (1998) and Bell and Williams (2001, 2005).

Harrison pioneered Brownian control models as a framework to find asymptotically optimal service policies for networks in heavy traffic. The framework was first proposed for multiclass queueing networks in Harrison (1988), and later was extended for stochastic processing networks in Harrison (2000). In his framework, a corresponding Brownian control problem of a stochastic processing network is first solved, and then the solution to the Brownian problem is used to construct service policies for the original stochastic processing network. Finally, these policies are shown to be asymptotically optimal for the stochastic processing network under a heavy traffic condition. A key step to solving the Brownian control problem is to have an equivalent workload formulation of the Brownian control problem as explained in Harrison and Van Mieghem (1997). The "workload process" of the Brownian control model corresponding to the stochastic processing network in this paper, as well as in Ata and Kumar (2005), Stolyar (2004), and Bell and Williams (2001, 2005), is one-dimensional. Thus, the Brownian control problem has a simple solution. Our maximum pressure policies, at least under a special linear holding cost structure, can be considered as another "interpretation" of the solution to the Brownian control problem, although this interpretation is not as direct as those in Bell and Williams (2005) and Ata and Kumar (2005). Our paper, together with these papers in the literature, demonstrates that the interpretation of the Brownian solution is not unique, proving the optimality of the interpreted policies can be difficult.

Maximum pressure type of policies were pioneered by Tassiluas and his coauthors under various names including back-pressure policies; see, for example, Tassiluas and Ephremides (1992, 1993), Tassiulas (1995) and Tassiulas and Bhattacharya (2000). The work of Tassiulas and Bhattacharya (2000) represents a significant advance in finding efficient operational policies for a wide class of networks, and is closely related to Dai and Lin (2005). Readers are referred to Dai and Lin (2005) for an explanation of the major differences of these two works. We note that, contrary to the description in Dai and Lin (2005), Tassiluas and Ephremides (1992, 1993) and Tassiulas (1995) do cover network models, not just one-pass systems. For a recent survey of these policies and their applications to wireless networks, see Georgiadis, Neely and Tassiulas (2006).

The remainder of the paper is organized as follows. In Section 1.1 we collect some of the notation used in this paper. In Section 2 we describe a class of stochastic processing networks, and introduce the maximum pressure service policies. We then define the workload process of a stochastic processing network in Section 3, where we also introduce the complete resource pooling condition. The main results of this paper are stated in Section 4. The proofs of the main theorems are outlined in Section 5. A key to the proofs of these theorems is a state space collapse result of the diffusion-scaled network processes under a maximum pressure policy. In Section 6 each fluid model solution under a maximum pressure policy is shown to exhibit a state space collapse. Section 7 applies Bramson's approach [Bramson (1998)] to prove the state space collapse of the diffusion-scaled network processes. The state space collapse result is converted into a heavy traffic limit theorem in Section 8. The limit theorem is used in Section 5 to complete the proofs of the main theorems. In Section 9 we discuss the $\varepsilon$-optimality of maximum pressure policies. A number of technical lemmas as well as Theorem 1 are proved in the Appendix A.

1.1. *Notation.* We use $\mathbb{R}^d$ to denote the $d$-dimensional Euclidean space. Vectors in $\mathbb{R}^d$ are envisioned as column vectors unless indicated otherwise. The transpose of a vector $v$ will be denoted as $v'$. For $v, w \in \mathbb{R}^d$, $v \cdot w$ denotes the dot product, and $v \times w$ denotes the vector $(v_1 w_1, \dots, v_d w_d)'$. The max norm in $\mathbb{R}^d$ is denoted as $|\cdot|$, and for a matrix $A$, we use $A$ to denote the maximum absolute value among all components. The Euclidean norm $\|\cdot\|$ in $\mathbb{R}^d$ is defined by $\|v\| = \sqrt{v \cdot v}$. For $r_1, r_2 \in \mathbb{R}$, we use $r_1 \vee r_2$ and $r_1 \wedge r_2$ to denote the maximum and minimum of $r_1$ and $r_2$, respectively.

We use $\mathbb{D}^d[0,\infty)$ to denote the set of functions $f : [0,\infty) \mapsto \mathbb{R}^d$ that are right continuous on $[0,\infty)$ having left limits in $(0,\infty)$. For $f \in \mathbb{D}^d[0,\infty)$, we let

$$\|f\|_t = \sup_{0 \le s \le t} |f(s)|.$$

We endow the function space $\mathbb{D}^d[0,\infty)$ with the usual Skorohod $J_1$-topology [Ethier and Kurtz (1986)]. A sequence of functions $\{f_r\} \subset \mathbb{D}^d[0,\infty)$ is said to converge to an $f \in \mathbb{D}^d[0,\infty)$ uniformly on compact (u.o.c.) sets, denoted as $f_r(\cdot) \to f(\cdot)$, if for each $t \geq 0$, $\lim_{r\to\infty} \|f_r - f\|_t = 0$. For a sequence of stochastic processes $\{X^r\}$ taking values in $\mathbb{D}^d[0,\infty)$, we use $X^r \Rightarrow X$ to denote that $X^r$ converges to $X$ in distribution.

## 2. Stochastic processing networks.

In this section we describe a general stochastic processing network proposed by Dai and Lin (2005). We follow the notation of Dai and Lin (2005). The network is assumed to have $\mathbf{I}+\mathbf{1}$ buffers, $\mathbf{J}$ activities and $\mathbf{K}$ processors. Buffers, activities and processors are indexed by $i = 0,\dots,\mathbf{I}$, $j = 1,\dots,\mathbf{J}$ and $k = 1,\dots,\mathbf{K}$, respectively. For notational convenience, we define $\mathcal{I} = \{1,\dots,\mathbf{I}\}$ the set of buffers excluding buffer 0, $\mathcal{J} = \{1,\dots,\mathbf{J}\}$ the set of activities and $\mathcal{K} = \{1,\dots,\mathbf{K}\}$ the set of processors. Each buffer, with infinite capacity, holds jobs or materials that await service. Buffer 0 is a special one that is used to model the outside world, where an infinite number of jobs await. Each activity can simultaneously process jobs from a set of buffers. It may require simultaneous possession of multiple processors to be active. Jobs departing from a buffer will go next to other buffers with certain probabilities that depend on the current activity taken.

### 2.1. *Resource consumption.*

Each activity needs one or more processors available to be active. For activity $j$, $A_{kj} = 1$, if activity $j$ requires processor $k$, and $A_{kj} = 0$ otherwise. The $\mathbf{K}\times\mathbf{J}$ matrix $A = (A_{kj})$ is the resource consumption matrix. Each activity may be allowed to process jobs in multiple buffers simultaneously. For activity $j$, we use the indicator function $B_{ji}$ to record whether buffer $i$ can be processed by activity $j$. ($B_{ji} = 1$ if activity $j$ processes buffer $i$ jobs.) The set of buffers $i$ with $B_{ji} = 1$ is said to be the *constituency* of activity $j$. It is denoted by $\mathcal{B}_j$. The constituency is assumed to be nonempty for each activity $j \in \mathcal{J}$, and may contain more than one buffer. When a processing requirement of an activity is met, a job departs from each one of the constituent buffers. For each activity $j$, we use $u_j(\ell)/\mu_j$ to denote the $\ell$th activity $j$ processing requirement, where $u_j = \{u_j(\ell), \ell \geq 1\}$ is an i.i.d. sequence of random variables, defined on some probability space $(\Omega, \mathcal{F}, \mathbb{P})$, and $\mu_j$ is a strictly positive real number. We set $\sigma_j^2 = \mathrm{var}(u_j(1))$, and assume that $\sigma_j < \infty$ and $u_j$ is *unitized*, that is, $\mathbb{E}[u_j(1)] = 1$, where $\mathbb{E}$ is the expectation operator associated with the probability measure $\mathbb{P}$. It follows that $1/\mu_j$ and $\sigma_j$ are the mean and coefficient of variation, respectively, for the processing times of activity $j$.

An activity $j$ is said to be an *input activity* if it processes jobs only from buffer 0, that is, $\mathcal{B}_j = \{0\}$. An activity $j$ is said to be a *service activity* if it does not process any job from buffer 0, that is, $0 \notin \mathcal{B}_j$. We assume that each activity is either an input activity or a service activity. We further assume that each processor processes either input activities only or service activities only. A processor that only processes input activities is called an *input processor*, and a processor that only processes service activities is called a *service processor*. The input processors process jobs from buffer 0 (outside) and generate the arrivals for the network. We denote $\mathcal{J}_I$ to be the set of input activities, $\mathcal{J}_S$ the set of service activities, $\mathcal{K}_I$ the set of input processors, and $\mathcal{K}_S$ the set of service processors.

2.2. *Routing.* Buffer $i$ jobs, after being processed by activity $j$, will go next to other buffers or leave the system. Let $e_0$ be the $\mathbf{I}$-dimensional vector of all 0's, and for $i \in \mathcal{I}$, $e_i$ is the $\mathbf{I}$-dimensional vector with the $i$th component 1 and other components 0. For each activity $j \in \mathcal{J}$ and each constituent buffer $i \in \mathcal{B}_j$, we use an $\mathbf{I}$-dimensional binary random vector $\phi_i^j(\ell) = (\phi_{ii'}^j(\ell), i' \in \mathcal{I})$ to denote the *routing vector* of the $\ell$th buffer $i$ job processed by activity $j$, where $\phi_i^j(\ell) = e_{i'}$ if the $\ell$th buffer $i$ job processed by activity $j$ goes next to buffer $i'$, and $\phi_i^j(\ell) = e_0$ if the job leaves the system. We assume that the sequence $\phi_i^j = \{\phi_i^j(\ell), \ell \geq 1\}$ is i.i.d., defined on the same probability space $(\Omega, \mathcal{F}, \mathbb{P})$, for each activity $j \in \mathcal{J}$ and $i \in \mathcal{B}_j$. Set $P_{ii'}^j = \mathbb{E}[\phi_{ii'}^j(1)]$. Then $P_{ii'}^j$ is the probability that a buffer $i$ job processed by activity $j$ will go next to buffer $i'$.

For each $j \in \mathcal{J}, i \in \mathcal{B}_j$, the cumulative routing process is defined by the sum

$$\Phi_i^j(\ell) = \sum_{n=1}^{\ell} \phi_i^j(n),$$

and $\Phi_{ii'}^j(\ell)$ denotes the number of jobs that will go next to buffer $i'$ among the first $\ell$ buffer $i$ jobs that are processed by activity $j$.

The sequences

$$(u_j, \phi_i^j : i \in \mathcal{B}_j, j \in \mathcal{J})$$

are said to be the *primitive increments* of the network. We assume that they are mutually independent and all are independent of the initial state of the network.

2.3. *Resource allocations.* Because multiple activities may require usage of the same processor, not all activities can be simultaneously undertaken at a 100% level. Unless stated otherwise, we assume that each processor's service capacity is infinitely divisible, and processor-splitting of a processor's service capacity is realizable. We use a nonnegative variable $a_j$ to denote the level at which processing activity $j$ is undertaken. When $a_j = 1$, activity $j$ is employed at a 100% level. When $a_j = 0$, activity $j$ is not employed. Suppose that the engagement level of activity $j$ is $a_j$, with $0 \le a_j \le 1$. The processing requirement of an activity $j$ job is depleted at rate $a_j$. (The job finishes processing when its processing requirement reaches 0.) The activity consumes $a_j A_{kj}$ fraction of processor $k$'s service capacity per unit time. The remaining service capacity, $1 - a_j A_{kj}$, can be used for other activities.

We use $a = (a_j) \in \mathbb{R}_+^{\mathbf{J}}$ to denote the corresponding $\mathbf{J}$-dimensional allocation (column) vector, where $\mathbb{R}_+$ denotes the set of nonnegative real numbers. Since each processor $k$ can decrease processing requirements at the rate of at most 1 per unit of time, we have

$$\sum_{j \in \mathcal{J}} A_{kj} a_j \le 1 \qquad \text{for each processor } k. \tag{2.1}$$

In vector form, $Aa \le e$, where $e$ is the $\mathbf{K}$-dimensional vector of ones. We assume that there is at least one input activity and that the input processors never idle. Namely,

$$\sum_{j \in \mathcal{J}} A_{kj} a_j = 1 \qquad \text{for each input processor } k. \tag{2.2}$$

We use $\mathcal{A}$ to denote the set of allocations $a \in \mathbb{R}_+^{\mathbf{J}}$ that satisfy (2.1) and (2.2).

Each $a \in \mathcal{A}$ represents an allowable allocation of the processors working on various activities. We note that $\mathcal{A}$ is bounded and convex. Let $\mathcal{E} = \{a^1, \ldots, a^{\mathbf{E}}\}$ be the set of extreme points of $\mathcal{A}$, where the total number $\mathbf{E}$ of extreme points is finite.

REMARK. Harrison (2002) does not have the concept of input processor and input activity for his stochastic processing networks. The input processors and activities in our network model allow us to capture the dynamic routing decisions for external arrivals. (Note that in our model the dynamic routing decision of a job is made not at the time of its arrival, but at the arrival time of the previous job.) However, if we restrict ourselves to the networks without these routing decisions, Harrison's network models are broader than ours; Harrison describes his models through first-order network data only, leaving the underlying stochastic primitives and system dynamics unspecified. In this sense, our network models are a special class of Harrison's, which include some stochastic processing networks where the service requirements for jobs from different buffers can be different even if they are processed simultaneously by a single activity. The latter networks are not covered in this paper.

2.4. *Service policies.* Each job in a buffer is assumed to be processed by one activity in its entire stay at the buffer. A processing of an activity can be preempted. In this case, each in-service job is "frozen" by the activity.

The next time the activity is made active again, the processing is resumed from where it was left off. In addition to the availability of processors, a (nonpreempted) activity can be made active only when each constituent buffer has jobs that are not in service or frozen. We assume that within each buffer jobs are queued in the order of their arrivals to the buffer, and *head-of-line* policy is used. When a (nonpreempted) activity becomes active with a given engagement level, the leading job in each buffer that is not in service or frozen is processed. If multiple activities are actively working on a buffer, there are multiple jobs in the buffer that are in service. For an allocation $a$, if there is an activity $j$ with $a_j > 0$ that cannot be made active, the allocation is infeasible. At any given time $t$, we use $\mathcal{A}(t)$ to denote the set of allocations that are *feasible* at that time. A service policy specifies which allocation is being undertaken at each time $t \geq 0$, and we use $\pi = \{\pi(t) : t \geq 0\}$ to denote such a policy. Under the policy $\pi$, allocation $\pi(t) \in \mathcal{A}(t)$ will be employed at time $t$.

Dai and Lin (2005) propose a family of service policies called *maximum pressure policies* that are throughput optimal for a large class of stochastic processing networks. To describe these policies for our network, for each buffer $i = 1, \ldots, \mathbf{I}$ and each activity $j = 1, \ldots, \mathbf{J}$, we define

$$R_{ij} = \mu_j \left( B_{ji} - \sum_{i' \in \mathcal{B}_j} P^j_{i'i} \right). \tag{2.3}$$

The $\mathbf{I} \times \mathbf{J}$ matrix $R = (R_{ij})$ is called the input-output matrix in Harrison (2002). One interprets $R_{ij}$ as the average amount of buffer $i$ material consumed per unit of activity $j$, with a negative value being interpreted to mean that activity $j$ is a net producer of material in buffer $i$. Define $\mathcal{E}(t) = \mathcal{E} \cap \mathcal{A}(t)$ to be the set of feasible extreme allocations at time $t$. Denote $Z = \{Z(t), t \geq 0\}$ to be the buffer level process with $Z_i(t)$ being the buffer level of buffer $i$ at time $t$, including those in service or "frozen" at time $t$. Now we are ready to define maximum pressure service policies for our stochastic processing network. Each maximum pressure policy is associated with a vector $\alpha \in \mathbb{R}^{\mathbf{I}}$ with $\alpha_i > 0$ for each $i \in \mathcal{I}$. Recall that for two vectors $x, y \in \mathbb{R}^d$, $x \times y$ denotes the vector $(x_1 y_1, \ldots, x_d y_d)'$.

Definition 1 (Maximum pressure policies). Given a vector $\alpha \in \mathbb{R}^{\mathbf{I}}$ with $\alpha_i > 0$ for $i \in \mathcal{I}$, a service policy is said to be a *maximum pressure policy associated with parameter* $\alpha$ if at each time $t$, the network chooses an allocation $a^* \in \arg\max_{a \in \mathcal{E}(t)} p_\alpha(a, Z(t))$, where $p_\alpha(a, Z(t))$ is called *the network pressure with parameter* $\alpha$ *under allocation* $a$ *and buffer level* $Z(t)$ and is defined as

$$p_\alpha(a, Z(t)) = (\alpha \times Z(t)) \cdot Ra. \tag{2.4}$$

When more than one allocation attains the maximum pressure, a tie-breaking rule is used. Our results are not affected by how ties are broken. However, for concreteness, one can order the extreme allocation set $\mathcal{E}$, and always choose the maximum-pressure allocation with the smallest index. In general, a maximum pressure policy can be an idling policy; namely, some processors may idle even if there are jobs that they can process.

In Section 4 we are going to show that every maximum pressure policy asymptotically minimizes a certain *workload process*, and when the stochastic processing network incurs a certain quadratic holding cost rate, there is a maximum pressure policy that is asymptotically optimal.

REMARK. Dai and Lin (2005) associate two parameters, $\alpha$ and $\beta$, with each maximum pressure policy and define the network pressure as $p_{\alpha,\beta}(a, Z(t)) = (\alpha \times Z(t) - \beta) \cdot Ra$. It turns out that the second parameter $\beta$ disappears in both fluid and diffusion limits and will have no impact on the asymptotic performance. Thus, for notational convenience, we set in this paper $\beta = 0$ and associate each maximum pressure policy with one parameter $\alpha$.

**3. Workload process and complete resource pooling.** We define the workload process through a linear program (LP) called the static planning problem and its dual problem. For a stochastic processing network with input-output matrix $R$ and capacity consumption matrix $A$, the static planning problem is defined as follows: choose a scalar $\rho$ and a $\mathbf{J}$-dimensional column vector $x$ so as to

$$\text{minimize} \quad \rho \tag{3.1}$$

$$\text{subject to} \quad Rx = 0, \tag{3.2}$$

$$\sum_{j \in \mathcal{J}} A_{kj} x_j = 1 \qquad \text{for each input processor } k, \tag{3.3}$$

$$\sum_{j \in \mathcal{J}} A_{kj} x_j \le \rho \qquad \text{for each service processor } k, \tag{3.4}$$

$$x \ge 0. \tag{3.5}$$

For each optimal solution $(\rho, x)$ to (3.1)–(3.5), the vector $x$ is said to be a *processing plan* for the stochastic processing network, where component $x_j$ is interpreted as the long-run fraction of time that activity $j$ is undertaken. Since one of the constraints in (3.4) must be binding for a service processor, $\rho$ is interpreted as the long-run utilization of the busiest service processor under the processing plan. With this interpretation, the left-hand side of (3.2) is interpreted as the long-run *net flow rates* from the buffers. Equality (3.2) demands that, for each buffer, the long-run input rate to the buffer is equal to the long-run output rate from the buffer. Equality (3.3) ensures that input processors never idle, while inequality (3.4) requires that each service processor's utilization not exceed that of the busiest service processor. The objective is to minimize the utilization of the busiest service processor. For future references, the optimal objective value $\rho$ is said to be the *traffic intensity* of the stochastic processing network.

The dual problem of the static planning problem is the following: choose an **I**-dimensional vector $y$ and a **K**-dimensional vector $z$ so as to

$$\text{maximize} \quad \sum_{k \in \mathcal{K}_I} z_k, \tag{3.6}$$

$$\text{subject to} \quad \sum_{i \in \mathcal{I}} y_i R_{ij} \le - \sum_{k \in \mathcal{K}_I} z_k A_{kj} \quad \text{for each input activity } j \tag{3.7}$$

$$\sum_{i \in \mathcal{I}} y_i R_{ij} \le \sum_{k \in \mathcal{K}_S} z_k A_{kj} \quad \text{for each service activity } j, \tag{3.8}$$

$$\sum_{k \in \mathcal{K}_S} z_k = 1, \tag{3.9}$$

$$z_k \ge 0 \quad \text{for each service processor } k. \tag{3.10}$$

Recall that $\mathcal{K}_I$ is the set of input processors, and $\mathcal{K}_S$ is the set of service processors. Each pair $(y, z)$ that satisfies (3.7)–(3.10) is said to be a *resource pool*. Component $y_i$ is interpreted as the work dedicated to a unit of buffer $i$ job by the resource pool, and $z_k$ is interpreted as the relative capacity of processor $k$, measured in fractions of the service capacity of the resource pool; for each input processor $k$, the relative capacity $z_k$ is the amount of work generated by input processor $k$ per unit of time. Equality (3.9) ensures that the service capacity of the resource pool equals the sum of service capacities of all service processors. Constraint (3.8) demands that no service activity can accomplish more work than the capacity it consumes. Recall that $-R_{ij}$ is the rate at which input activity $j$ generates buffer $i$ jobs. For each input activity $j$, constraint (3.7), which can be written as $\sum_{i \in \mathcal{I}} y_i(-R_{ij}) \ge \sum_{k \in \mathcal{K}_I} z_k A_{kj}$, ensures that the work dedicated to each unit of the activity is no less than that which it generates. The objective is to maximize $\sum_{k \in \mathcal{K}_I} z_k$, which is the total amount of work generated from outside by the input processors per unit of time. A service processor $k$ is said to be in the resource pool $(y, z)$ if $z_k > 0$. Constraint (3.9) is an equality instead of an inequality because we do not pose a nonnegativity constraint on the variable $\rho$ in the primal LP; the nonnegativity of $\rho$ is guaranteed from constraint (3.4) and the nonnegativity of $A$ and $x$.

A *bottleneck pool* is defined to be an optimal solution $(y^*, z^*)$ to the dual LP (3.6)–(3.10). Let $(\rho^*, x^*)$ be an optimal solution to the primal LP, the static planning problem (3.1)–(3.5). From the basic duality theory, $\sum_j A_{kj} x_j^* = \rho^*$ for any service processor $k$ with $z_k^* > 0$. It says that all service processors in the bottleneck pool $(y^*, z^*)$ are the busiest servers under any optimal processing plan $x^*$.

For a bottleneck pool $(y^*, z^*)$, let $W(t) = y^* \cdot Z(t)$ for $t \geq 0$. Then, $W(t)$ represents the average total work of this bottleneck pool embodied in all jobs that are present at time $t$ in the stochastic processing network. The process $W = \{W(t), t \geq 0\}$ is called the *workload process* of this bottleneck pool. Although the workload process of a nonbottleneck resource pool $(y, z)$ can also be defined by $y \cdot Z(t)$, we will focus on the workload processes of bottleneck pools because bottleneck pools become significantly more important in heavy traffic. In general, the bottleneck pool is not unique. However, we assume all the stochastic processing networks considered in this paper have a unique bottleneck pool; namely, they satisfy the following *complete resource pooling condition*.

DEFINITION 2 (Complete resource pooling condition). A stochastic processing network is said to satisfy the complete resource pooling condition if the corresponding dual static planning problem (3.6)–(3.10) has a nonnegative, unique optimal solution $(y^*, z^*)$.

For a processing network that satisfies the complete resource pooling condition, we define the *bottleneck workload process*, or simply the workload process, of the stochastic processing network to be the workload process of its unique bottleneck pool.

The (bottleneck) workload process defined here is different from the workload process defined in Harrison and Van Mieghem (1997). Their workload process is multi-dimensional, with some components corresponding to the nonbottleneck pools; it is defined in terms of what they call "reversible displacements." For the networks where their workload process has dimension one, these two definitions of the workload process are consistent.

REMARK. Under certain assumptions including a heavy traffic assumption that requires all servers in the network be critically loaded, Harrison (2000) proposes a "canonical" representation of the workload process for a stochastic processing network through a dual LP that is similar to (3.6)–(3.10). There, basic optimal solutions to the dual LP were chosen as rows of the workload matrix which was used to define the workload process. Without his heavy traffic assumption, his "canonical" choice of workload matrix would exclude those nonbottleneck servers. In this case, it is not clear how to define a "canonical" representation of the workload process to include those nonbottleneck servers. Although for some network examples like multiclass queueing networks we can define the workload matrix such that its rows are the basic solutions to the dual LP, more analysis is required for general stochastic processing networks.

**4. Asymptotic optimality and the main results.** The behavior of the buffer level process and the workload process for a stochastic processing network under any policy is complex. In particular, deriving closed form expressions for performance measures involving these processes is not possible in general. Therefore, we perform an asymptotic analysis for stochastic processing networks operating under maximum pressure policies. Our asymptotic region is when the network is in heavy traffic; that is, the offered traffic load is approximately equal to the system capacity. Formally, we consider a sequence of stochastic processing networks indexed by $r = 1, 2, \ldots$; as $r \to \infty$, the traffic intensity $\rho^r$ of the $r$th network goes to one. We assume that these networks all have the same network topology and primitive increments. In other words, the matrices $A$ and $B$, and the sequences $(u_j, \phi_i^j : j \in \mathcal{J}, i \in \mathcal{B}_j)$ do not vary with $r$. However, we allow the processing rates to change with $r$, and use $\mu_j^r$ to denote the processing rate of activity $j$ in the $r$th network. Thus, the traffic intensity $\rho^r$ for the $r$th network is the optimal objective value of the static planning problem (3.1)–(3.5) with the input-output matrix $R^r = (R_{ij}^r)$ given by $R_{ij}^r = \mu_j^r(B_{ji} - \sum_{i' \in \mathcal{B}_j} P_{i'i}^j)$. We assume the following *heavy traffic assumption* throughout this paper.

ASSUMPTION 1 (Heavy traffic assumption). There exists a constant $\mu_j > 0$ for each activity $j \in \mathcal{J}$ such that, as $r \to \infty$,

$$\mu_j^r \to \mu_j, \tag{4.1}$$

and, setting $R = (R_{ij})$ as in (2.3) with $\mu_j$ being the limit values in (4.1), the static planning problem (3.1)–(3.5) with parameter $(R, A)$ has a unique optimal solution $(\rho^*, x^*)$ with $\rho^* = 1$. Furthermore, as $r \to \infty$,

$$r(\rho^r - 1) \to \theta \tag{4.2}$$

for some constant $\theta$.

We define the *limit network* of the network sequence to be the network that has the same network topology and primitive increments as networks in the sequence, and has processing rates equal to the limit values $\mu_j$, given in (4.1). Assumption 1 basically means that in the limit network there exists a unique processing plan $x^*$ that can avoid inventory buildups over time, and the busiest service processor is fully utilized under this processing plan. Condition (4.2) requires that the networks' traffic intensities approach to 1 at rate $r^{-1}$ or faster.

The heavy traffic assumption is now quite standard in the literature; see, for example, Bramson and Dai (2001), Chen and Zhang (2000), Dai and Kurtz (1995) and Williams (1996, 1998a) for heavy traffic analysis of queueing networks and Ata and Kumar (2005), Harrison (2000), Harrison and López (1999), Harrison and Van Mieghem (1997) and Williams (2000) for heavy traffic analysis of stochastic processing networks. However, heavy traffic assumptions in the literature usually assume that, in addition to Assumption 1, *all* service processors are fully utilized. The latter assumption, together with a complete resource pooling assumption to be introduced later in this section, rules out some common networks such as multiclass queueing networks. In our heavy traffic assumption, only the busiest service processor is required to be critically loaded, and some other service processors are allowed to be under-utilized.

The optimal processing plan $x^*$ given in Assumption 1 is referred to as the *nominal processing plan*. We use $T_j(t)$ to denote the cumulative amount of activity $j$ processing time in $[0,t]$ for the limit network; let $T(t)$ be the corresponding $\mathbf{J}$-dimensional vector. Then $T(t)/t$ is the average activity levels over the time span $[0,t]$. To avoid a linear buildup of jobs over time in the limit network, the long-run average rate (or activity level) that activity $j$ is undertaken needs to equal $x_j^*$, that is,

$$\lim_{t\to\infty} T(t)/t = x^* \qquad \text{almost surely.} \tag{4.3}$$

There should be no linear buildup of jobs under a reasonably "good" policy. A policy is said to be efficient for the limit network if (4.3) holds for the network operating under the policy. Since we consider a sequence of networks, we would like to define an analogous notion of a "good" or efficient policy for the sequence. One can imagine that under a reasonably "good" policy, when $r$ is large, the average activity levels over long time spans must be very close to the nominal processing plan $x^*$. To be specific, we define the notion of *asymptotic efficiency* as follows. Let $T_j^r(t)$ be the cumulative amount of activity $j$ processing time in $[0,t]$ for the $r$th network and $T^r(t)$ be the corresponding $\mathbf{J}$-dimensional vector.

Definition 3 (Asymptotic efficiency). Consider a sequence of stochastic processing networks indexed by $r = 1, 2, \dots,$ where Assumption 1 holds. A policy $\pi$ is said to be *asymptotically efficient* if and only if under policy $\pi$, with probability 1, for each $t \geq 0$,

$$T^r(r^2 t)/r^2 \to x^* t \qquad \text{as } r \to \infty. \tag{4.4}$$

Equation (4.4) basically says that, under an asymptotically efficient policy, the average activity levels over a time span of order $r^2$ are very close to the nominal processing plan, so that no linear buildup of jobs will occur over the time span of this order.

Asymptotic efficiency is closely related to the *throughput optimality* as defined in Dai and Lin (2005). Fluid models have been used to prove throughput optimality of stochastic processing networks operating under a policy.

Similarly, the fluid model corresponding to the limit network can be used to prove the asymptotic efficiency of a policy for the sequence of networks that satisfies Assumption 1. In particular, one can prove that a policy $\pi$ is asymptotically efficient if the fluid model of the limit network operating under $\pi$ is weakly stable. To introduce the first result of this paper, we first define the extreme-allocation-available (EAA) condition for a stochastic processing network.

DEFINITION 4 (EAA condition). A stochastic processing network is said to satisfy the *EAA condition* if, for any vector $q \in \mathbb{R}^{\mathbf{I}}_{+}$, there exists an extreme allocation $a^* \in \mathcal{E}$ such that $Ra^* \cdot q = \max_{a\in\mathcal{E}} Ra \cdot q$, and for each buffer $i$ with $\sum_j a^*_j B_{ji} > 0$, buffer level $q_i$ is positive.

Readers are referred to Section 6 of Dai and Lin (2005) for a detailed discussion on the EAA condition; there a class of networks is shown to satisfy the EAA condition and a network example is shown not to satisfy the EAA condition.

THEOREM 1. *Consider a sequence of stochastic processing networks that satisfies Assumption 1. If the limit network satisfies the EAA condition, then for any $\alpha \in \mathbb{R}^{\mathbf{I}}$ with $\alpha > 0$, the maximum pressure policy with parameter $\alpha$ is asymptotically efficient.*

Using a fluid model, Dai and Lin (2005) prove in their Theorem 2 that, for a stochastic processing network satisfying the *EAA condition*, a maximum pressure policy is throughput optimal. The proof of Theorem 1 uses a fluid model that will be introduced in Section 6 and is almost identical to that of Theorem 2 in Dai and Lin (2005). It will be outlined in Appendix B.

Asymptotic efficiency helps to identify reasonably "good" policies, but it is not very discriminating. We would like to demonstrate a certain sense of optimality for maximum pressure policies in terms of secondary performance measures. For this, we will introduce two notions of *asymptotic optimality*. Before giving their definitions, we make the following assumption on the sequence of networks.

ASSUMPTION 2 (Complete resource pooling). All networks in the sequence and the limit network satisfy the complete resource pooling condition defined in Section 3. Namely, the dual static planning problem (3.6)–(3.10) of the $r$th network has a nonnegative, unique optimal solution $(y^r, z^r)$, and the dual static planning problem of the limit network also has a nonnegative, unique optimal solution $(y^*, z^*)$.

LEMMA 1. *Assume Assumptions 1 and 2. Then $(y^r, z^r) \to (y^*, z^*)$ as $r \to \infty$.*

The proof of Lemma 1 will be provided in Appendix B.

REMARK. In Assumption 2 we assume all networks in the sequence satisfy the complete resource pooling condition so that the workload processes $W^r$ can be uniquely defined as in Section 3 by the first order network data $(R^r, A^r)$. This assumption can be removed if one defines the workload process of a network with multiple bottleneck pools to be the workload process of an arbitrarily chosen but prespecified bottleneck pool (with $y^r$ being any given optimal solution to the dual problem). On the other hand, the complete resource pooling condition for the limit network is crucial for our results to hold.

The first notion of asymptotic optimality is in terms of the workload process introduced in Section 3. Under Assumption 2, we can define the one-dimensional workload process of the $r$th network as

$$W^r(t) = y^r \cdot Z^r(t), \tag{4.5}$$

where $Z^r = \{Z^r(t), t \geq 0\}$ is the buffer level process of the $r$th network.

Define the diffusion-scaled workload and buffer level processes of the $r$th network $\widehat{W}^r = \{\widehat{W}^r(t), t \geq 0\}$ and $\widehat{Z}^r = \{\widehat{Z}^r(t), t \geq 0\}$ via

$$\widehat{W}^r(t) = W^r(r^2 t)/r, \qquad \widehat{Z}^r(t) = Z^r(r^2 t)/r.$$

Clearly, $\widehat{W}^r(t) = y^r \cdot \widehat{Z}^r(t)$ for $t \geq 0$.

DEFINITION 5. Consider a sequence of stochastic processing networks indexed by $r$. An asymptotically efficient policy $\pi^*$ is said to be *asymptotically optimal for workload* if for any $t > 0, w > 0$, and any asymptotically efficient policy $\pi$,

$$\limsup_{r\to\infty} \mathbb{P}(\widehat{W}^r_{\pi^*}(t) > w) \leq \liminf_{r\to\infty} \mathbb{P}(\widehat{W}^r_{\pi}(t) > w), \tag{4.6}$$

where $\widehat{W}^r_{\pi^*}(\cdot)$ and $\widehat{W}^r_{\pi}(\cdot)$ are the diffusion-scaled workload processes under policies $\pi^*$ and $\pi$, respectively.

The second notion of asymptotic optimality is in terms of a quadratic holding cost structure for the sequence of stochastic processing networks. Let $h = (h_i : i \in \mathcal{I})$ be a constant vector with $h_i > 0$ for $i \in \mathcal{I}$. For the $r$th network, the diffusion-scaled holding cost rate at time $t$ is

$$\widehat{H}^r(t) = \sum_{i\in\mathcal{I}} h_i(\widehat{Z}^r_i(t))^2. \tag{4.7}$$

DEFINITION 6. Consider a sequence of stochastic processing networks indexed by $r$. An asymptotically efficient policy $\pi^*$ is said to be *asymptotically optimal for the quadratic holding cost* if for any $t>0, \eta>0$, and any asymptotically efficient policy $\pi$,

$$\limsup_{r\to\infty} \mathbb{P}(\widehat{H}^r_{\pi^*}(t)>\eta) \le \liminf_{r\to\infty} \mathbb{P}(\widehat{H}^r_{\pi}(t)>\eta), \tag{4.8}$$

where $\widehat{H}^r_{\pi^*}(t)$ and $\widehat{H}^r_{\pi}(t)$ are the diffusion-scaled total holding cost rates at time $t$ under policies $\pi^*$ and $\pi$, respectively.

To state our main theorems, we make two more assumptions on the sequence of networks. One is a moment assumption on the unitized service times $u_j(\ell)$ and the other is an assumption on the initial buffer level processes.

ASSUMPTION 3. There exists an $\varepsilon_u>0$ such that, for all $j$,

$$\mathbb{E}[(u_j(1))^{2+\varepsilon_u}]<\infty. \tag{4.9}$$

ASSUMPTION 4 (Initial condition). As $r\to\infty$,

$$\widehat{Z}^r(0)=Z^r(0)/r\to 0 \qquad \text{in probability.} \tag{4.10}$$

Assumption 3 requires that the unitized service times have finite $2+\varepsilon_u$ moments. It was used by Ata and Kumar (2005) and it is stronger than some standard regularity assumptions such as in Bramson (1998). Assumption 3 will be used in Section 7 to prove a state space collapse result for stochastic processing networks operating under maximum pressure policies. Assumption 4 holds if the initial buffer levels of the networks are stochastically bounded, namely,

$$\lim_{\tau\to\infty}\limsup_{r\to\infty}\mathbb{P}(|Z^r(0)|>\tau)=0.$$

Clearly, Assumption 4 implies that, as $r\to\infty$,

$$\widehat{W}^r(0)\to 0 \qquad \text{in probability.} \tag{4.11}$$

THEOREM 2. *Consider a sequence of stochastic processing networks. Assume Assumptions 1–4 and that the limit network satisfies the EAA condition. Each maximum pressure policy is asymptotically optimal for workload.*

THEOREM 3. *Consider a sequence of stochastic processing networks where Assumptions 1–4 hold and the limit network satisfies the EAA condition. The maximum pressure policy with parameter $\alpha=h$ is asymptotically optimal for the quadratic holding cost in (4.7).*

Theorem 2 says that, at every time $t$, asymptotically, the diffusion-scaled workload under any maximum pressure policy $\pi^*$ is dominated by that under any other asymptotically efficient policy $\pi$ in the sense of stochastic ordering. Theorem 3 says that, given a quadratic holding cost rate structure with coefficient vector $h$, the quadratic holding cost rate at every time $t$ under the maximum pressure policy with parameter $\alpha = h$ is asymptotically dominated by that under any other efficient policy in the sense of stochastic ordering.

The proofs of Theorems 2 and 3 will be outlined in Section 5. Throughout this paper, we shall assume Assumptions 1–4 and that the limit network satisfies the EAA condition.

When the objective is to minimize a linear holding cost, our maximum pressure policies do not achieve the asymptotic optimality. However, some maximum pressure policies are asymptotically near optimal. In Section 9 we introduce the notion of $\varepsilon$-optimality and identify a set of maximum pressure policies that are *asymptotically $\varepsilon$-optimal* in terms of minimizing the linear holding cost.

**5. An outline of the proofs.** This section outlines the proofs of our main theorems, Theorems 2 and 3. We first derive an asymptotic lower bound on the workload processes under asymptotically efficient policies. Then we state a heavy traffic limit result in Theorem 4, which implies that this asymptotic lower bound is achieved by the workload process under any maximum pressure policy. The heavy traffic limit theorem also implies the optimality of a certain maximum pressure policy under the quadratic holding cost structure in Theorem 3. At the end we outline a proof for the heavy traffic limit theorem. The key to the proof is a state space collapse result to be stated in Theorem 5.

We first derive an asymptotic lower bound on the workload processes under asymptotically efficient policies. That is, we search for a process $W^*$ such that, under any asymptotically efficient policy,

$$\liminf_{r\to\infty} \mathbb{P}(\widehat{W}^r(t) > w) \geq \mathbb{P}(W^*(t) > w) \qquad \text{for all } t \text{ and } w.$$

As before, we assume that the sequence of stochastic processing networks satisfies Assumption 2. Hence, the one dimensional workload process $W^r$ is well defined by (4.5) for each $r$.

We begin the analysis by defining a process $Y^r = \{Y^r(t), t \geq 0\}$ for the $r$th network via

$$Y^r(t) = (1 - \rho^r)t - y^r \cdot R^r T^r(t). \tag{5.1}$$

Since $\rho^r$ is interpreted as the traffic intensity of the bottleneck pool, for each $t \geq 0$, $\rho^r t$ is interpreted as the average total work contributed to the bottleneck pool from the exogenous arrivals in $[0, t]$, and $(1 - \rho^r)t$ represents the average total work that could have been depleted by time $t$ if the bottleneck pool never idles. Because of the randomness of the processing times, the bottleneck pool will almost surely incur idle time over time, particularly when the system is not overloaded. Under a service policy and its corresponding activity level process $T^r$, the average total work that has been depleted by time $t$ is given by

$$y^r \cdot R^r T^r(t) = \sum_{j\in\mathcal{J}_S} T_j^r(t) \sum_{i\in\mathcal{I}} y_i^r R_{ij}^r - \sum_{j\in\mathcal{J}_I} T_j^r(t) \sum_{i\in\mathcal{I}} y_i^r(-R_{ij}^r).$$

Note that, as in Section 3, for each service activity $j \in \mathcal{J}_S$, $\sum_{i\in\mathcal{I}} y_i^r R_{ij}^r$ is the average work accomplished per unit of activity $j$, and that for each input activity $j \in \mathcal{J}_I$, $\sum_{i\in\mathcal{I}} y_i^r(-R_{ij}^r)$ is the average work generated per unit of activity $j$. Therefore, $Y^r(t)$ represents the deviation of the workload depletion in $[0,t]$ from that under the "best" policy. The following lemma says that this deviation does not decrease over time.

LEMMA 2. *Consider a sequence of networks satisfying Assumption 2. For each $r$ and each sample path, the process $Y^r$ defined in (5.1) is a nondecreasing function with $Y^r(0) = 0$.*

We leave the proof to Appendix B.

REMARK. Some special stochastic processing networks, such as multiclass queueing networks [Harrison (1988)] and unitary networks [Bramson and Williams (2003)], have no control on the input activities. Then, $T_j^r(t)$ is fixed for all $j \in \mathcal{J}_I$ under different policies, and $\sum_{j\in\mathcal{J}_I} T_j^r(t) \sum_{i\in\mathcal{I}} y_i^r(-R_{ij}^r) = \rho^r t$. For these networks, one gets

$$Y^r(t) = t - \sum_{j\in\mathcal{J}_S} T_j^r(t) \sum_{i\in\mathcal{I}} y_i^r R_{ij}^r,$$

and $Y^r(t)$ is interpreted as the cumulative idle time of the bottleneck pool by time $t$.

For each activity $j$, we define the process $S_j^r = \{S_j^r(t), t \geq 0\}$ by

$$S_j^r(t) = \max\left\{n : \sum_{\ell=1}^{n} u_j(\ell) \leq \mu_j^r t\right\}.$$

Under a head-of-line service policy, $S_j^r(t)$ is the number of activity $j$ processing completions in $t$ units of activity $j$ processing time for the $r$th network. Then we can describe the system dynamics of the $r$th network by the following equation:

$$\begin{aligned} Z_i^r(t) = Z_i^r(0) &+ \sum_{j\in\mathcal{J}} \sum_{i'\in\mathcal{B}_j} \Phi_{i'i}^j(S_j^r(T_j^r(t))) \\ &- \sum_{j\in\mathcal{J}} S_j^r(T_j^r(t)) B_{ji} \qquad \text{for each } t \geq 0 \text{ and } i \in \mathcal{I}. \end{aligned} \tag{5.2}$$

Since quantity $T_j^r(t)$ is the cumulative amount of activity $j$ processing time in $[0,t]$, $S_j^r(T_j^r(t))$ is the number of activity $j$ processings completed by time $t$, and $\sum_{j\in\mathcal{J}} S_j^r(T_j^r(t))B_{ji}$ is the total number of jobs that depart from buffer $i\in\mathcal{I}\cup\{0\}$ in $[0,t]$. For each activity $j$, $\sum_{i'\in\mathcal{B}_j}\Phi_{i'i}^j(S_j^r(T_j^r(t)))$ is the total number of jobs sent to buffer $i$ by activity $j$ from its constituent buffers by time $t$, so $\sum_{j\in\mathcal{J}}\sum_{i'\in\mathcal{B}_j}\Phi_{i'i}^j(S_j^r(T_j^r(t)))$ is the total number of jobs that go to buffer $i$ by time $t$. Equation (5.2) says that the number of jobs in buffer $i$ at time $t$ equals the initial number plus the number of arrivals minus the number of departures.

From (5.2), we can write the workload process $W^r = y^r\cdot Z^r$ as

$$W^r(t) = W^r(0) + \sum_{i\in\mathcal{I}} y_i^r \sum_{j\in\mathcal{J}}\left(\sum_{i'\in\mathcal{B}_j}\Phi_{i'i}^j(S_j^r(T_j^r(t))) - B_{ji}S_j^r(T_j^r(t))\right).$$

Let $X^r(t) = W^r(t) - Y^r(t)$ for each $t\ge 0$. Then

$$\begin{aligned}X^r(t) = W^r(0) &+ \sum_{i\in\mathcal{I}} y_i^r \sum_{j\in\mathcal{J}}\left(\sum_{i'\in\mathcal{B}_j}\Phi_{i'i}^j(S_j^r(T_j^r(t))) - B_{ji}S_j^r(T_j^r(t))\right)\\ &- (1-\rho^r)t + y^r\cdot R^rT^r(t).\end{aligned}$$

We define the following diffusion-scaled processes:

$$\hat{S}_j^r(t) = r^{-1}[S_j^r(r^2t) - \mu_j^r r^2 t] \qquad \text{for each } j\in\mathcal{J}, \tag{5.3}$$

$$\hat{\Phi}_{i'i}^{j,r}(t) = r^{-1}[\Phi_{i'i}^j(\lfloor r^2 t\rfloor) - P_{i'i}^j r^2 t] \qquad \text{for each } j\in\mathcal{J}, i'\in\mathcal{B}_j \text{ and } i\in\mathcal{I}, \tag{5.4}$$

$$\hat{X}^r(t) = r^{-1}X^r(r^2t), \tag{5.5}$$

$$\hat{Y}^r(t) = r^{-1}Y^r(r^2t). \tag{5.6}$$

In (5.4) $\lfloor t\rfloor$ denotes the greatest integer number less than or equal to the real number $t$.

Then the diffusion-scaled workload process $\widehat{W}^r$ can be written as a sum of two processes

$$\widehat{W}^r(t) = \hat{X}^r(t) + \hat{Y}^r(t), \tag{5.7}$$

and

$$\begin{aligned}\hat{X}^r(t) = \widehat{W}^r(0) &+ \sum_{i\in\mathcal{I}} y_i^r \sum_{j\in\mathcal{J}}\left(\sum_{i'\in\mathcal{B}_j}\hat{\Phi}_{i'i}^{j,r}(S_j^r(T_j^r(t)))\right.\\ &\left. + \left(\sum_{i'\in\mathcal{B}_j} P_{i'i}^j - B_{ji}\right)\hat{S}_j^r(T_j^r(t))\right)\\ &- r(1-\rho^r)t,\end{aligned} \tag{5.8}$$

where

$$\bar{\bar{T}}^r_j(t) = r^{-2}T^r_j(r^2t) \quad \text{and} \quad S^r_j(t) = r^{-2}S^r_j(r^2t).$$

The process $\widehat{X}^r$ depends on the policy only through the fluid-scaled process $T^r$. In fact, from Lemma 4.1 of Dai (1995) and (4.4), it follows that, under any asymptotically efficient policy, $T^r \to x^*(\cdot)$ u.o.c. almost surely, where $x^*(t) = x^*t$ and $x^*$ is the optimal solution to the static planning problem (3.1)–(3.5) of the limit network. As a consequence, $\widehat{X}^r$ converges in distribution to a one-dimensional Brownian motion that is independent of policies.

LEMMA 3. *Consider a sequence of stochastic processing networks operating under an asymptotically efficient policy. Assume Assumptions 1–4. Then $\widehat{X}^r \Rightarrow X^*$, where $X^*$ is a one-dimensional Brownian motion that starts from the origin, has drift parameter $\theta$ given in (4.2), and has variance parameter*

$$\sigma^2 = (y^*)'\left(\sum_{j\in\mathcal{J}} x^*_j\mu_j \sum_{i\in\mathcal{B}_j} \Upsilon^{j,i}\right)y^* + \sum_{i\in\mathcal{I}}\sum_{j\in\mathcal{J}}(y^*_i)^2R^2_{ij}x^*_j\mu_j\sigma^2_j$$

*with $\Upsilon^{j,i}, j\in\mathcal{J}, i\in\mathcal{B}_j$, defined by*

$$\Upsilon^{j,i}_{i_1,i_2} = \begin{cases} P^j_{i,i_1}(1-P^j_{i,i_2}), & \text{if } i_1 = i_2,\\ -P^j_{i,i_1}P^j_{i,i_2}, & \text{if } i_1 \neq i_2.\end{cases}$$

PROOF. First, Lemma 4.1 of Dai (1995) and (4.4) implies that $\bar{\bar{T}}^r(\cdot) \to x^*(\cdot)$ almost surely under any asymptotically efficient policy. Then, the result in the lemma follows from (5.8), the functional central limit theorem for renewal processes [cf. Iglehart and Whitt (1970)] the random time change theorem [cf. Billingsley (1999), (17.9)], and the continuous mapping theorem [cf. Billingsley (1999), Theorem 5.1]. Deriving the expression for $\sigma^2$ is straightforward but tedious; the derivation is outlined in Harrison (1988) for multiclass queueing networks, so we will not repeat it here. □

We define the one-dimensional reflection mapping $\psi:\mathbb{D}[0,\infty)\to\mathbb{D}[0,\infty)$ such that, for each $f\in\mathbb{D}[0,\infty)$ with $f(0)\geq 0$,

$$\psi(f)(t) = f(t) - \inf_{0\leq s\leq t}(f(s)\wedge 0).$$

Applying diffusion scaling to Lemma 2, we know that $\widehat{Y}^r(\cdot)$ is a nonnegative, nondecreasing function, so, from (5.7) and the well-known minimality of the solution of the one-dimensional Skorohod problem [cf. Bell and Williams (2001), Proposition B.1],

$$\widehat{W}^r(t) \geq \psi(\widehat{X}^r)(t) \qquad \text{for every } t \text{ and every sample path;}$$

namely, $\psi(\widehat{X}^r)(t)$ is a pathwise lower bound on $\widehat{W}^r$. It then follows that

$$\liminf_{r\to\infty} \mathbb{P}(\widehat{W}^r(t) > w) \ge \liminf_{r\to\infty} \mathbb{P}(\psi(\widehat{X}^r)(t) > w) \qquad \text{for all } t \text{ and } w.$$

Define

$$W^* = \psi(X^*). \tag{5.9}$$

Then $W^*$ is a one-dimensional reflecting Brownian motion associated with $X^*$. Because $\widehat{X}^r \Rightarrow X^*$, by the continuous mapping theorem, we have

$$\psi(\widehat{X}^r) \quad \Rightarrow \quad W^*.$$

Because $W^*(t)$ has continuous distribution for each $t$, we have

$$\lim_{r\to\infty} \mathbb{P}(\psi(\widehat{X}^r)(t) > w) = \mathbb{P}(W^*(t) > w).$$

Therefore,

$$\liminf_{r\to\infty} \mathbb{P}(\widehat{W}^r(t) > w) \ge \mathbb{P}(W^*(t) > w) \qquad \text{for each } t \text{ and } w. \tag{5.10}$$

So far, we have shown that $W^*$ is an asymptotic lower bound on the workload processes under asymptotically efficient policies. The following heavy traffic limit theorem ensures that the workload processes under any maximum pressure policy converge to $W^*$. This completes the proof of Theorem 2.

THEOREM 4. *Consider a sequence of stochastic processing networks operating under a maximum pressure policy with parameter $\alpha > 0$. Assume Assumptions 1–4 and that the limit network satisfies the EAA condition. Then*

$$(\widehat{W}^r, \widehat{Z}^r) \quad \Rightarrow \quad (W^*, Z^{*,\alpha}) \qquad \text{as } r \to \infty, \tag{5.11}$$

*where $W^*$ is given by (5.9) and*

$$Z^{*,\alpha} = \zeta^\alpha W^*, \tag{5.12}$$

*with $\zeta^\alpha$ being defined as*

$$\zeta_i^\alpha = \frac{y_i^*/\alpha_i}{\sum_{i'} (y_{i'}^*)^2/\alpha_{i'}}. \tag{5.13}$$

An outline for proving Theorem 4 will be presented at the end of this section. The full proof of Theorem 4 will be completed in Section 8. Theorem 4 is known as a *heavy traffic limit theorem*. It states that the scaled buffer level process and the scaled workload process jointly converge in heavy traffic as displayed in (5.11), and the limit processes exhibit a form of state space collapse: in diffusion limit, the **I**-dimensional buffer level process is a constant vector multiple of the one-dimensional workload process as displayed in (5.12). Because buffers at nonbottleneck pools may also contribute to the workload process (at the bottleneck pool), the diffusion-scaled buffer level processes in these buffers do not go to zero in heavy traffic limit, due to the *idling* nature of our maximum pressure policies. For example, consider a network of two stations in series, known as the tandem queue network. Suppose that the second station is critically loaded and the first station is underloaded. Under any maximum pressure policy, by Theorem 4, the diffusion limit for the buffer level process at the first station is a one-dimensional reflecting Brownian motion. This result is in sharp contrast to the ones in Chen and Mandelbaum (1991), in which they prove buffers at nonbottleneck stations in a single-class network disappear in heavy traffic diffusion limit under the first-come–first-serve *nonidling* policy.

With Theorem 4, we now give the proof of Theorem 3.

Proof of Theorem 3. We first derive an asymptotic lower bound on the quadratic holding cost rate process under any asymptotically efficient policy. Consider the following quadratic optimization problem:

$$g^r(w) = \min_{q\geq 0} \sum_{i\in\mathcal{I}} h_i q_i^2 \tag{5.14}$$

$$\text{s.t.} \quad \sum_{i\in\mathcal{I}} y_i^r q_i = w. \tag{5.15}$$

The problem (5.14)–(5.15) can be solved easily. The optimal solution is given by

$$q_i^r = (y_i^r/h_i) \Big/ \left(\sum_{i'\in\mathcal{I}} (y_{i'}^r)^2/h_{i'}\right) \qquad \text{for } i\in\mathcal{I},$$

and the optimal objective value is given by

$$g^r(w) = w^2 \Big/ \left(\sum_{i'\in\mathcal{I}} (y_{i'}^r)^2/h_{i'}\right).$$

Then under any policy, for all $r$, we have

$$\widehat{H}^r(t) = \sum_{i\in\mathcal{I}} h_i(\widehat{Z}_i^r(t))^2 \geq g^r(\widehat{W}^r(t)) = (\widehat{W}^r(t))^2 \Big/ \left(\sum_{i'\in\mathcal{I}} (y_{i'}^r)^2/h_{i'}\right).$$

For any asymptotically efficient policy $\pi$,

$$\begin{aligned}\liminf_{r\to\infty} \mathbb{P}(\widehat{H}_\pi^r(t) > \eta) &\geq \liminf_{r\to\infty} \mathbb{P}\left((\widehat{W}_\pi^r(t))^2 \Big/ \left(\sum_{i'\in\mathcal{I}} (y_{i'}^r)^2/h_{i'}\right) > \eta\right) \\ &\geq \mathbb{P}(g^*(W^*(t)) > \eta),\end{aligned}$$

where $g^*(w) = w^2/(\sum_{i'\in\mathcal{I}}(y^*_{i'})^2/h_{i'})$. The second inequality follows from (5.10) and the fact that $y^r \to y^*$ as $r \to \infty$.

Now consider the maximum pressure policy $\pi^*$ having the parameter $\alpha = h$. From Theorem 4, we have

$$\widehat{H}^r_{\pi^*} \quad\Rightarrow\quad \sum_{i\in\mathcal{I}} h_i(Z_i^{*,h})^2 \qquad \text{as } r \to \infty.$$

Because

$$\sum_i h_i(Z_i^{*,h}(t))^2 = \sum_i h_i(\zeta_i^h W^*(t))^2 = g^*(W^*(t)),$$

we have

$$\lim_{r\to\infty} \mathbb{P}(\widehat{H}^r_{\pi^*}(t) > \eta) = \mathbb{P}(g^*(W^*(t)) > \eta),$$

thus proving inequality (4.8). □

The key to the proof of Theorem 4 is the following state space collapse result.

THEOREM 5 (State space collapse). *Consider a sequence of stochastic processing networks operating under the maximum pressure policy with parameter $\alpha$. Assume Assumptions 1–4 and that the limit network satisfies the EAA condition. Then, for each $T \geq 0$, as $r \to \infty$,*

$$\|\widehat{Z}^r(\cdot) - \zeta^\alpha \widehat{W}^r(\cdot)\|_T \to 0 \qquad \textit{in probability.} \tag{5.16}$$

Recall that $\|\cdot\|_T$ is the uniform norm over $[0,T]$. [The readers should not confuse the symbols $T$ and $T(\cdot)$ with one another. We will always include "$(\cdot)$" when dealing with the cumulative activity level process $T(\cdot)$.] Theorem 5 states a form of state space collapse for the diffusion-scaled network process: for large $r$, the $\mathbf{I}$-dimensional diffusion-scaled buffer level process is essentially a constant vector multiple of the one-dimensional workload process.

The proof of Theorem 5 is lengthy. To prove it, we generalize a framework of Bramson (1998) from the multiclass queueing network setting to the stochastic processing network setting. The framework consists of two steps that we will follow in the next two sections: First, in Section 6 we will show that any fluid model solution for the stochastic processing networks under a maximum pressure policy exhibits some type of state space collapse, which is stated in Theorem 6 in that section. Then, in Section 7 we will follow Bramson's approach to translate the state space collapse of the fluid model into a state space collapse result under diffusion scaling, and thus prove Theorem 5.

Once we have Theorem 5, we apply a perturbed Skorohod mapping theorem from Williams (1998b) to complete the proof of Theorem 4 in Section 8.

**6. The fluid model.** In this section we first introduce the fluid model of a sequence of stochastic processing networks operating under a maximum pressure policy. We then show that any fluid model solution under a maximum pressure policy exhibits a form of state space collapse.

The fluid model of a sequence of stochastic processing networks that satisfies Assumption 1 is defined by the following equations: for all $t \ge 0$,

$$\bar{Z}(t) = \bar{Z}(0) - R\bar{T}(t), \tag{6.1}$$

$$\bar{Z}(t) \ge 0, \tag{6.2}$$

$$\sum_{j \in \mathcal{J}} A_{kj}(\bar{T}_j(t) - \bar{T}_j(s)) = t - s \quad \text{for each } 0 \le s \le t \text{ and each input processor } k, \tag{6.3}$$

$$\sum_{j \in \mathcal{J}} A_{kj}(\bar{T}_j(t) - \bar{T}_j(s)) \le t - s \quad \text{for each } 0 \le s \le t \text{ and each processor } k, \tag{6.4}$$

$$\bar{T} \text{ is nondecreasing and } \bar{T}(0) = 0. \tag{6.5}$$

Equations (6.1)–(6.5) define the fluid model under *any* given service policy. Any quantity $(\bar{Z}, \bar{T})$ that satisfies (6.1)–(6.5) is a *fluid model solution* to the fluid model that operates under a general service policy. Following its stochastic processing network counterparts, each fluid model solution $(\bar{Z}, \bar{T})$ has the following interpretations: $\bar{Z}_j(t)$ the fluid level in buffer $i$ at time $t$ and $\bar{T}_j(t)$ the cumulative amount of activity $j$ processing time in $[0, t]$.

Under a specific service policy, there are additional fluid model equations. For a given parameter $\alpha > 0$, we are to specify the fluid model equation associated with the maximum pressure policy with parameter $\alpha$. To motivate the fluid model equation, we note that for each fluid model solution $(\bar{Z}, \bar{T})$, it follows from equations (6.3)–(6.4) that $\bar{T}$, and hence $\bar{Z}$, is Lipschitz continuous. Thus, the solution is absolutely continuous and has derivatives almost everywhere with respect to the Lebesgue measure on $[0, \infty)$. A time $t > 0$ is said to be a regular point of the fluid model solution if the solution is differentiable at time $t$. For a function $f : \mathbb{R}_+ \to \mathbb{R}^d$, where $d$ is some positive integer, we use $\dot{f}(t)$ to denote the derivative of $f$ at time $t$ when the derivative exists. From (6.3)–(6.4), one has $\dot{\bar{T}}(t) \in \mathcal{A}$ at each regular time $t$. Thus,

$$R\dot{\bar{T}}(t) \cdot (\alpha \times \bar{Z}(t)) \le \max_{a \in \mathcal{A}} Ra \cdot (\alpha \times \bar{Z}(t)) = \max_{a \in \mathcal{E}} Ra \cdot (\alpha \times \bar{Z}(t)).$$

Assume that the limit network satisfies the EAA condition. The fluid model equation associated with the maximum pressure policy with parameter $\alpha$ takes the following form:

$$R\dot{\bar{T}}(t)\cdot(\alpha\times\bar{Z}(t))=\max_{a\in\mathcal{E}}Ra\cdot(\alpha\times\bar{Z}(t)) \tag{6.6}$$

for each regular time $t$.

Each fluid model equation will be justified through a fluid limit procedure. Two types of fluid limits are considered in this paper. One will be introduced in Section 7 and the other in Appendix B. They both satisfy the fluid model equations (6.1)–(6.6). Equation (6.6) says that, under the maximum pressure policy with parameter $\alpha$, the instantaneous activity allocation $\dot{\bar{T}}(t)$ in the fluid model maximizes the network pressure, $Ra\cdot(\alpha\times\bar{Z}(t))$, at each regular time $t$. Any fluid model solution that satisfies fluid model equations (6.1)–(6.6) is called a *fluid model solution under the maximum pressure policy* with parameter $\alpha$.

THEOREM 6. *Let $(\bar{Z},\bar{T})$ be a solution to the fluid model equations (6.1)–(6.6), where $\alpha>0$, and $(R,A)$ is the first order network data for the limit network that satisfies Assumptions 1 and 2. Suppose $|\bar{Z}(0)|\le 1$. Then there exists some finite $\tau_0>0$, which depends on just $\alpha$, $\mathbf{I}$, $R$ and $A$, such that,*

$$|\bar{Z}(t)-\zeta^\alpha\bar{W}(t)|=0 \qquad \text{for all } t\ge\tau_0, \tag{6.7}$$

*where $\bar{W}=y^*\cdot\bar{Z}$ is the workload process of the fluid model and $\zeta^\alpha$ is given by (5.13). Furthermore, if*

$$|\bar{Z}(\tau_1)-\zeta^\alpha\bar{W}(\tau_1)|=0 \qquad \text{for some } \tau_1\ge 0,$$

*then*

$$|\bar{Z}(t)-\zeta^\alpha\bar{W}(t)|=0 \qquad \text{for all } t\ge\tau_1.$$

Theorem 6 says that the fluid model under the maximum pressure policy exhibits a form of state space collapse: after some finite time $\tau_0$, the $\mathbf{I}$-dimensional buffer level process $\bar{Z}$ equals a constant vector multiple of the one-dimensional workload process $\bar{W}$; if this happens at time $\tau_1$, it happens at all times after $\tau_1$. In particular, if $\bar{Z}(0)=\zeta^\alpha\bar{W}(0)$, then $\bar{Z}(t)=\zeta^\alpha\bar{W}(t)$ for all $t\ge 0$.

The rest of this section is devoted to the proof of Theorem 6. We first define

$$\bar{Z}^*(t)=\zeta^\alpha\bar{W}(t),$$

and we shall prove $\bar{Z}(t)-\bar{Z}^*(t)=0$ for $t$ large enough. The following lemma will be used in the proof.

LEMMA 4. *For any* $t \geq 0$,

$$(\bar{Z}(t) - \bar{Z}^*(t)) \cdot y^* = 0, \tag{6.8}$$

$$(\bar{Z}(t) - \bar{Z}^*(t)) \cdot (\alpha \times \bar{Z}^*(t)) = 0, \tag{6.9}$$

*and for each regular time* $t$,

$$(\bar{Z}(t) - \bar{Z}^*(t)) \cdot (\alpha \times \dot{\bar{Z}}^*(t)) = 0. \tag{6.10}$$

PROOF. Equality (6.8) follows because

$$\bar{Z}^*(t) \cdot y^* = \bar{W}(t)\zeta^\alpha \cdot y^* = \bar{W}(t) = \bar{Z}(t) \cdot y^*.$$

Because $\bar{Z}^*(t) = \zeta^\alpha \bar{W}(t)$,

$$\alpha \times \bar{Z}^*(t) = \bar{W}(t) y^* \Big/ \left( \sum_{i \in \mathcal{I}} (y_i^*)^2 / \alpha_i \right). \tag{6.11}$$

Then equality (6.9) follows from (6.8).

Finally, equality (6.10) follows immediately from (6.9). □

Define $V = \{Ra : a \in \mathcal{A}\}$. Recall that $\mathcal{A}$ is the set of all possible allocations and the vector $Ra$ is the average rate at which material is consumed from all buffers under allocation $a$, so $V$ is the set of all possible flow rates out of buffers in the limit network. It is obvious that $V$ is a polytope containing the origin because $Rx^* = 0$, where, as before, $x^*$ is the optimal solution to the static planning problem of the limit network. Furthermore, Corollary A.1 in Appendix A, applied to the limit network, implies that

$$\max_{v \in V} y^* \cdot v = 0. \tag{6.12}$$

Since $y^* \cdot Ra$ is the rate at which the workload is reduced under allocation $a$, therefore, (6.12) says that no feasible allocation can reduce the system workload for the limit network, which is not surprising because the limit network is in heavy traffic.

Define a $(\mathbf{I} - 1)$-dimensional hyperplane

$$V^o = \{v \in \mathbb{R}^I : y^* \cdot v = 0\}. \tag{6.13}$$

The hyperplane $V^o$ is a supporting hyperplane of $V$. Let $V^* = V \cap V^o$. Then from (6.12) and (6.13), $V^* = \arg\max_{v \in V} y^* \cdot v$; namely, $V^*$ is the set of all possible flow rates that maximize $y^* \cdot v$. Since the origin is in both $V$ and $V^o$, $V^*$ is not empty. Moreover, the following lemma says that $V^*$ contains a certain $(\mathbf{I} - 1)$-dimensional neighborhood of the origin.

LEMMA 5. *There exists some* $\delta > 0$ *such that* $\{v \in V^o : \|v\| \leq \delta\} \subset V$.

We provide a proof of Lemma 5 in Appendix B. A key to the proof is the complete resource pooling assumption for the limit network; namely, the dual static planning problem (3.6)–(3.10) has a unique optimal solution. The uniqueness of $y^*$ ensures that the origin lies in the *relative* interior of one of the facets of $V$.

PROOF OF THEOREM 6. Let $(\bar{Z}, \bar{T})$ be a solution to the fluid model equations (6.1)–(6.6). We consider the Lyapunov function

$$f(t) = (\alpha \times (\bar{Z}(t) - \bar{Z}^*(t))) \cdot (\bar{Z}(t) - \bar{Z}^*(t)).$$

Let $v(t) = R\dot{\bar{T}}(t)$; for each buffer $i$, $v_i(t)$ can be interpreted as the net flow rate out of buffer $i$ at time $t$ (*total departure rate minus total arrival rate*). Then it follows from (6.1) that $\dot{\bar{Z}}(t) = -v(t)$, and we have

$$\begin{aligned} \dot{f}(t) &= 2(\alpha \times (\bar{Z}(t) - \bar{Z}^*(t))\Big) \cdot (\dot{\bar{Z}}(t) - \dot{\bar{Z}}^*(t)) \\ &= 2(\alpha \times (\bar{Z}(t) - \bar{Z}^*(t))) \cdot (-v(t)). \end{aligned} \tag{6.14}$$

The second equality in (6.14) follows from (6.10) in Lemma 4.

Because $\dot{\bar{T}}(t) \in \mathcal{A}$, we have $v(t) \in V$. Thus, $y^* \cdot v(t) \le 0$. This, together with (6.11), implies that $(\alpha \times \bar{Z}^*(t)) \cdot v(t) \le 0$. Furthermore, the fluid model equation (6.6) implies that

$$(\alpha \times \bar{Z}(t)) \cdot v(t) = \max_{v \in V} (\alpha \times \bar{Z}(t)) \cdot v.$$

Therefore, the last term in (6.14) is bounded from above as follows:

$$\begin{aligned} 2(\alpha \times (\bar{Z}(t) - \bar{Z}^*(t))) \cdot (-v(t)) &\le -2(\alpha \times \bar{Z}(t)) \cdot v(t) \\ &= -2 \max_{v \in V} (\alpha \times \bar{Z}(t)) \cdot v. \end{aligned} \tag{6.15}$$

Since $V^* \subset V$, we have

$$\begin{aligned} \max_{v \in V} (\alpha \times \bar{Z}(t)) \cdot v &\ge \max_{v \in V^*} (\alpha \times \bar{Z}(t)) \cdot v \\ &= \max_{v \in V^*} (\alpha \times (\bar{Z}(t) - \bar{Z}^*(t))) \cdot v. \end{aligned} \tag{6.16}$$

The second equality in (6.16) holds because of (6.11) and the fact that $y^* \cdot v = 0$ for all $v \in V^*$.

If $f(t) > 0$, let

$$v^* = \frac{\delta(\bar{Z}(t) - \bar{Z}^*(t))}{\|\bar{Z}(t) - \bar{Z}^*(t)\|}.$$

Then $\|v^*\| = \delta$ and $y^* \cdot v^* = 0$. The latter fact follows from (6.8). Pick $\delta$ as in Lemma 5 and it follows from Lemma 5 that $v^* \in V^*$. Therefore,

$$\begin{aligned}\max_{v\in V^*}(\alpha \times (\bar{Z}(t) - \bar{Z}^*(t))) \cdot v &\geq (\alpha \times (\bar{Z}(t) - \bar{Z}^*(t))) \cdot v^* \\ &= \delta f(t)/\|\bar{Z}(t) - \bar{Z}^*(t)\|.\end{aligned} \tag{6.17}$$

Combining (6.14)–(6.17), we have

$$\dot{f}(t) \leq -2\delta f(t)/\|\bar{Z}(t) - \bar{Z}^*(t)\| \leq -2\delta \sqrt{\min_{i\in\mathcal{I}} \alpha_i}\sqrt{f(t)}. \tag{6.18}$$

Therefore,
$f(t) = 0$ for $t \geq \sqrt{f(0)}/(\delta\sqrt{\min_{i\in\mathcal{I}} \alpha_i})$. Set $\tau_0 = \sqrt{\mathbf{I}\max_{i\in\mathcal{I}\alpha_i}}/(\delta\sqrt{\min_{i\in\mathcal{I}} \alpha_i})$. Then $f(t) = 0$ for $t \geq \tau_0$, because

$$\begin{aligned} f(0) &= \alpha \times (\bar{Z}(0) - \bar{Z}^*(0)) \cdot (\bar{Z}(0) - \bar{Z}^*(0)) = \alpha \times (\bar{Z}(0) - \bar{Z}^*(0)) \cdot \bar{Z}(0) \\ &\leq \max_{i\in\mathcal{I}} \alpha_i \|\bar{Z}(0)\|^2 \leq \mathbf{I} \max_{i\in\mathcal{I}} \alpha_i. \end{aligned}$$

Here $\tau_0$ depends only on $\alpha$, $R$, $A$ and $\mathbf{I}$, because the set $V$ is completely determined by $R$ and $A$ and so is $\delta$.

Equation (6.18) also implies that for any $\tau_1$, if $f(\tau_1) = 0$, then $f(t) = 0$ for all $t \geq \tau_1$. $\square$

**7. State space collapse.** In this section we translate the state space collapse result of the fluid model into a state space collapse result under the diffusion scaling, thus proving Theorem 5. We apply Bramson's approach [Bramson (1998)] to prove that, for each $T \geq 0$, as $r \to \infty$,

$$\|\widehat{Z}^r(\cdot) - \zeta^\alpha \widehat{W}^r(\cdot)\|_T \to 0 \qquad \text{in probability.}$$

In Bramson's approach, the following fluid scaling plays an important role in connecting Theorem 6 with Theorem 5: For each $r = 1, 2, \ldots,$ and $m = 0, 1, \ldots,$

$$S_j^{r,m}(t) = \frac{1}{\xi_{r,m}}(S_j^r(rm + \xi_{r,m}t) - S_j^r(rm)) \qquad \text{for each } j \in \mathcal{J},$$

$$\Phi_i^{j,r,m}(t) = \frac{1}{\xi_{r,m}}(\Phi_i^j(S_j^r(rm) + \lfloor \xi_{r,m}t \rfloor) - \Phi_i^j(S_j^r(rm))) \qquad \text{for each } j \in \mathcal{J}, i \in \mathcal{B}_j,$$

$$T_j^{r,m}(t) = \frac{1}{\xi_{r,m}}(T_j^r(rm + \xi_{r,m}t) - T_j^r(rm)) \qquad \text{for each } j \in \mathcal{J},$$

$$Z_i^{r,m}(t) = \frac{1}{\xi_{r,m}} Z_i^r(rm + \xi_{r,m}t) \qquad \text{for each } i \in \mathcal{I},$$

where $\xi_{r,m} = |Z^r(rm)| \vee r$. Recall that $\lfloor t \rfloor$ denotes the integer part of $t$.

Here scaling the processes by $\xi_{r,m}$ ensures $|Z^{r,m}(0)| \le 1$, which is needed for compactness reasons. Using index $(r,m)$ allows the time scale to expand; we will examine the processes over $[0,L]$ for $m = 0,1,\ldots,\lceil rT \rceil - 1$, where $L > 0$ and $T > 0$ are fixed. Thus, the diffusion-scaled time $[0,T]$ is covered by $\lceil rT \rceil$ fluid scaled time pieces, each with length $L$. Here $\lceil t \rceil$ denotes the smallest integer greater than or equal to $t$. For future references, we assume that $L > \tau_0 \vee 1$, where $\tau_0$ is the time when state space collapse occurs in the corresponding fluid model as determined in Theorem 6.

We outline the proof of Theorem 5 as follows. First, in Proposition 1 below, we give a probability estimate on the upper bound of the fluctuation of the stochastic network processes $\mathbb{X}^{r,m}(\cdot) = (T^{r,m}(\cdot), Z^{r,m}(\cdot))$. The estimates on the service processes $S_j^{r,m}(\cdot)$ and the routing processes $\Phi_i^{j,r,m}(\cdot)$ are also given. From Proposition 1, a so-called "good" set $\mathcal{G}^r$ of sample paths can be defined, where the processes $\mathbb{X}^{r,m}$ perform nicely for $r$ large enough. On this "good" set, for large enough $r$, the processes $\mathbb{X}^{r,m}$ can be uniformly approximated by so-called *Lipschitz cluster points*. These cluster points will be shown to be fluid model solutions under the maximum pressure policy. Since the state space collapse result holds for the fluid model under the maximum pressure policy by Theorem 6, the network processes $\mathbb{X}^{r,m}$ asymptotically have the state space collapse. The latter result will be translated into the state space collapse result for diffusion-scaled processes $\widehat{\mathbb{X}}^r$ as $r$ approaches infinity.

Note that in Theorem 6 the state space collapse of the fluid model does not happen instantaneously after time 0 if the initial state does not exhibit a state space collapse. The fluid-scaled processes $\mathbb{X}^{r,m}$ start from time $rm$ in the original network processes. Hence, for $m \ge 1$, $\mathbb{X}^{r,m}$ do not automatically have the state space collapse at the initial point. Such collapse can only be expected to occur in the interval $[\tau_0, L]$. However, for $m = 0$, the state space collapse happens at time 0 because of the initial condition (4.10), therefore, the collapse remains in the whole interval $[0,L]$. For this reason, we separate the proof into two parts according to the two intervals in the diffusion-scaled time: $[0, L\xi_{r,0}/r^2]$ and $[\tau_0\xi_{r,0}/r^2, T]$. Note that the two intervals overlap because of our assumption on the choice of $L$; if $\tau_0\xi_{r,0}/r^2 > T$, the second interval disappears and the first interval covers the whole interval $[0,T]$.

Propositions 2–5 develop a version of state space collapse for the interval $[\tau_0\xi_{r,0}/r^2, T]$. Proposition 2 shows that, on $\mathcal{G}^r$, the scaled processes $\mathbb{X}^{r,m}(\cdot)$ are uniformly close to Lipschitz cluster points for large $r$. Proposition 3 shows that the above cluster points are solutions to the fluid model equations. In Proposition 4 we use Propositions 2 and 3 and Theorem 6 to prove the state space collapse of the fluid scaled processes $\mathbb{X}^{r,m}(\cdot)$ in the time interval $[\tau_0, L]$ on the "good" set $\mathcal{G}^r$. The result in Proposition 4 is then translated into diffusion scaling in Proposition 5 to give a version of state space collapse for the diffusion-scaled process $\hat{\mathbb{X}}^r(t)$ in the interval $[\tau_0\xi_{r,0}/r^2, T]$.

In Propositions 6–8, we show that the state space collapse occurs for the diffusion-scaled process $\hat{\mathbb{X}}^r(t)$ in the interval $[0, L\xi_{r,0}/r^2]$. The basic idea is the same as the one described in the preceding paragraph except that now we only consider the scaled processes with $m = 0$. The corresponding network processes start from time 0, and by assuming the state space collapse happens at time 0 in Assumption 4, we have a stronger result for these types of processes: the state space collapse holds during the whole time interval $[0, L]$ instead of just in $[\tau_0, L]$. In fact, the scaled processes $\mathbb{X}^{r,0}(\cdot)$ are proved to be uniformly close to some cluster points for which the state space collapse starts at time 0. These facts are stated in Propositions 6 and 7. In Proposition 8 we summarize the state space collapse result for the fluid-scaled process $\mathbb{X}^{r,0}(\cdot)$ on $[0, L]$ and translate it into a result for the diffusion-scaled process in $[0, L\xi_{r,0}/r^2]$.

The results to be obtained in Propositions 5 and 8 are actually a *multiplicative state space collapse*, as called in Bramson (1998). To obtain the state space collapse result stated in Theorem 5, we will prove that $\xi_{r,m}/r$ are stochastically bounded at the end of this section.

7.1. *Probability estimates.* In this section we give probability estimates on the service processes $S_j^{r,m}(\cdot)$, the routing processes $\Phi_i^{j,r,m}(\cdot)$ and the upper bound of the fluctuation of the stochastic network processes $\mathbb{X}^{r,m}(\cdot)$.

PROPOSITION 1. *Consider a sequence of stochastic processing networks where the moment assumption, Assumption 3, is assumed. Fix $\epsilon > 0$, $L$ and $T$. Then, for large enough $r$,*

$$\mathbb{P}\left(\max_{m<rT} \|S_j^{r,m}(T_j^{r,m}(t)) - \mu_j^r T_j^{r,m}(t)\|_L > \epsilon\right) \le \epsilon \quad \text{for each } j \in \mathcal{J}, \tag{7.1}$$

$$\mathbb{P}\left(\max_{m<rT} \|\Phi_i^{j,r,m}(S_j^{r,m}(T_j^{r,m}(t))) - P_i^j \mu_j^r(T_j^{r,m}(t))\|_L > \epsilon\right) \le \epsilon \quad \text{for each } j \in \mathcal{J} \text{ and } i \in \mathcal{B}_j, \tag{7.2}$$

$$\mathbb{P}\left(\sup_{0\le t_1\le t_2\le L} |\mathbb{X}^{r,m}(t_2) - \mathbb{X}^{r,m}(t_1)| > N|t_2 - t_1| + \epsilon \text{ for some } m < rT\right) \le \epsilon, \tag{7.3}$$

*where $N = \mathbf{J}(1 + |R|)$.*

Recall that $|R| = \max_{ij} R_{ij}$. The following lemma is essential to the proof of Proposition 1.

LEMMA 6. *Assume that the moment assumption, Assumption 3, holds. Then for given $T$, and each $j \in \mathcal{J}$,*

$$u_j^{r,T,\max}/r \to 0 \qquad \text{as } r \to \infty \text{ with probability } 1, \tag{7.4}$$

*where* $u_j^{r,T,\max} = \max\{u_j(\ell) : 1 \le \ell \le S_j^r(r^2T) + 1\}$. *Furthermore, for any given* $\epsilon > 0$,

$$\mathbb{P}(\|\Phi_i^j(\ell) - P_i^j \ell\|_n \ge \epsilon n) \le \epsilon/n \qquad \text{for each } j \in \mathcal{J} \text{ and } i \in \mathcal{B}_j, \text{ and large enough } n, \tag{7.5}$$

*and for large enough* $t$,

$$\mathbb{P}(\|S_j^r(\tau) - \mu_j^r \tau\|_t \ge \epsilon t) \le \epsilon/t \qquad \text{for all } j \in \mathcal{J}, \text{and all } r. \tag{7.6}$$

We delay the proof of Lemma 6 to Appendix B, and we now prove Proposition 1.

PROOF OF PROPOSITION 1. The proof here essentially follows the same reasoning as in Propositions 5.1 and 5.2 of Bramson (1998). We first investigate the processes with index $m$, and then multiply the error bounds by the number of processes in each case, $\lceil rT \rceil$. We start with (7.1).

Fix a $j \in \mathcal{J}$, and pick any $\epsilon_1 \in (0,1]$. Let $\tau$ be the time that the first activity $j$ service completion occurs after time $rm$. We consider two cases. First suppose that $\tau > rm + \xi_{r,m} T_j^{r,m}(t)$. Then $\xi_{r,m} T_j^{r,m}(t) < \tau - rm \le u_j^{r,T,\max}/\mu_j^r$ and $S_j^r(rm + \xi_{r,m} T_j^{r,m}(t)) = S_j^r(rm)$. Thus,

$$\|S_j^r(rm + \xi_{r,m} T_j^{r,m}(t)) - S_j^r(rm) - \mu_j^r \xi_{r,m} T_j^{r,m}(t)\|_L \le u_j^{r,T,\max}.$$

Or if $\tau \le rm + \xi_{r,m} T_j^{r,m}(t)$, then by restarting the process $S_j^r$ at time $\tau$, we have, by (7.6), that for large enough $r$,

$$\begin{aligned}
&\mathbb{P}(\|S_j^r(rm + \xi_{r,m} T_j^{r,m}(t)) - S_j^r(\tau) - \mu_j^r (rm + \xi_{r,m} T_j^{r,m}(t) - \tau)\|_L \ge \epsilon_1 L \xi_{r,m}) \\
&\quad \le \epsilon_1/(Lr);
\end{aligned}$$

we use the fact that $T_j^{r,m}(t) \le t$ and $\xi_{r,m} \ge r$. Because $S_j^r(\tau) = S_j^r(rm) + 1$ and $\tau - rm \le u_j^{r,T,\max}/\mu_j^r$, we have

$$\begin{aligned}
&\mathbb{P}(\|S_j^r(rm + \xi_{r,m} T_j^{r,m}(t)) \\
&\qquad - S_j^r(rm) - \mu_j^r \xi_{r,m} T_j^{r,m}(t)\|_L \ge |1 - u_j^{r,T,\max}| + \epsilon_1 L \xi_{r,m}) \\
&\qquad \le \epsilon_1/(Lr).
\end{aligned}$$

In both cases we have, for large enough $r$,

$$\begin{aligned}
&\mathbb{P}(\|S_j^{r,m}(T_j^{r,m}(t)) - \mu_j^r T_j^{r,m}(t)\|_L \ge (1 + u_j^{r,T,\max})/\xi_{r,m} + \epsilon_1 L) \\
&\qquad \le \epsilon_1/(Lr).
\end{aligned} \tag{7.7}$$

From (7.4), we have, for large enough $r$,

$$\mathbb{P}(u_j^{r,T,\max}/r \geq \epsilon_1) \leq \epsilon_1. \tag{7.8}$$

Let $\mathcal{M}^r$ denote the complement of the events in (7.8). Then, for large enough $r$,

$$\mathbb{P}(\mathcal{M}^r) \geq 1 - \epsilon_1. \tag{7.9}$$

Then we have, for large enough $r$,

$$\begin{aligned}
&\mathbb{P}[(\|S_j^{r,m}(T_j^{r,m}(t)) - \mu_j^r T_j^{r,m}(t)\|_L \geq (L+2)\epsilon_1) \cap \mathcal{M}^r] \\
&\quad \leq \mathbb{P}[(\|S_j^{r,m}(T_j^{r,m}(t)) - \mu_j^r T_j^{r,m}(t)\|_L \\
&\qquad \geq (1 + u_j^{r,T,\max})/\xi_{r,m} + \epsilon_1 L) \cap \mathcal{M}^r].
\end{aligned} \tag{7.10}$$

In (7.10) we let $r$ be large enough so that $1/\xi_{r,m} \leq \epsilon_1$ and, hence, $(1 + u_j^{r,T,\max})/\xi_{r,m} \leq 2\epsilon_1$ for all events in $\mathcal{M}^r$.

It then follows from (7.7) and (7.10) that

$$\mathbb{P}[(\|S_j^{r,m}(T_j^{r,m}(t)) - \mu_j^r T_j^{r,m}(t)\|_L \geq (L+2)\epsilon_1) \cap \mathcal{M}^r] \leq \epsilon_1/(Lr). \tag{7.11}$$

Inequality (7.11) holds for $m = 1, \dots, \lceil rT \rceil$. It then follows that

$$\begin{aligned}
&\mathbb{P}\left[\left(\max_{m \leq rT} \|S_j^{r,m}(T_j^{r,m}(t)) - \mu_j^r T_j^{r,m}(t)\|_L \geq (L+2)\epsilon_1\right) \cap \mathcal{M}^r\right] \\
&\quad \leq \lceil rT \rceil \epsilon_1/(Lr) \leq \epsilon_1 T.
\end{aligned}$$

Then we have

$$\begin{aligned}
&\mathbb{P}\left(\max_{m \leq rT} \|S_j^{r,m}(T_j^{r,m}(t)) - \mu_j^r T_j^{r,m}(t)\|_L \geq (L+2)\epsilon_1\right) \\
&\quad \leq \epsilon_1 T + (1 - \mathbb{P}(\mathcal{M}^r)) \leq (T+1)\epsilon_1.
\end{aligned}$$

Since $\epsilon_1$ can be chosen arbitrarily, we let $\epsilon_1 = 1 \wedge \epsilon/((L+2) \vee (T+1))$. Then (7.1) follows.

Let $\mu^{\max}$ be an upper bound on $|\mu^r|$. Then from (7.11), for large enough $r$,

$$\mathbb{P}[(S_j^{r,m}(T_j^{r,m}(L)) \geq (\mu^{\max} L + L + 2)) \cap \mathcal{M}^r] \leq \epsilon_1/(Lr).$$

From (7.5), replacing $n$ by $(\mu^{\max} L + L + 2)\xi_{r,m}$, we have for each $i \in \mathcal{B}_j$,

$$\begin{aligned}
&\mathbb{P}[(\|\Phi_i^{j,r,m}(S_j^{r,m}(T_j^{r,m}(t))) - P_i^j S_j^{r,m}(T_j^{r,m}(t))\|_L \\
&\qquad > \epsilon_1(\mu^{\max} L + L + 2)) \cap \mathcal{M}^r] \\
&\quad \leq \epsilon_1/((\mu^{\max} L + L + 2)r) + \epsilon_1/(Lr).
\end{aligned} \tag{7.12}$$

Let $P^{\max} = \max_{j\in\mathcal{J}, i\in\mathcal{B}_j} P_i^j$. Then from (7.11), we have

$$\mathbb{P}[(\|P_i^j S_j^{r,m}(T_j^{r,m}(t)) - P_i^j \mu_j^r T_j^{r,m}(t)\|_L \geq P^{\max}(L+2)\epsilon_1) \cap \mathcal{M}^r] \\ \leq \epsilon_1/(Lr). \tag{7.13}$$

It follows from (7.12) and (7.13) that

$$\begin{aligned}\mathbb{P}[(\|\Phi_i^{j,r,m}(S_j^{r,m}(T_j^{r,m}(t))) - P_i^j \mu_j^r(T_j^{r,m}(t))\|_L \\ > \epsilon_1(P^{\max} + \mu^{\max} + 1)(L+2)) \cap \mathcal{M}^r] \\ \leq 5\epsilon_1/(2Lr).\end{aligned}$$

This inequality holds for $m = 1, \ldots, \lceil rT \rceil$, so setting $\epsilon_1 = 1 \wedge \epsilon/[((P^{\max} + \mu^{\max} + 1)(L+2)) \vee (5T+1)]$, one gets

$$\mathbb{P}\left[\left(\max_{m\leq rT} \|\Phi_i^{j,r,m}(S_j^{r,m}(T_j^{r,m}(t))) - P_i^j \mu_j^r(T_j^{r,m}(t))\|_L > \epsilon\right) \cap \mathcal{M}^r\right] \leq 5T\epsilon_1.$$

Thus, (7.2) follows.

Now we are going to show (7.3). First, it is easy to see that, for each $j \in \mathcal{J}$ and each $r$,

$$T_j^r(t) - T_j^r(s) \leq t - s \qquad \text{for } 0 \leq s \leq t$$

along any sample path. Therefore, the bounds in (7.3) on components $T_j$ are obtained with $N = 1$. For components $Z_i, i \in \mathcal{I}$, scaling (5.2) and applying (7.1) and (7.2) gives that for any $\epsilon_2 > 0$ and large enough $r$,

$$\begin{aligned}\mathbb{P}\Bigg(\sup_{0\leq t_1\leq t_2\leq L} |Z_i^{r,m}(t_2) - Z_i^{r,m}(t_1)| \\ > \sum_{j\in\mathcal{J}} (|R_{ij}^r|)|T_j^{r,m}(t_2) - T_j^{r,m}(t_1)| + 2\mathbf{J}(\mathbf{I}+2)\epsilon_2 \text{ for some } m \leq rT\Bigg) \\ \leq 2\mathbf{J}(\mathbf{I}+2)\epsilon_2.\end{aligned} \tag{7.14}$$

Since $R^r \to R$ as $r \to \infty$, we can choose $r$ large enough so that $|R^r| \leq |R| + 1$. Then, setting $\epsilon_2 = \epsilon/(2\mathbf{J}(\mathbf{I}+2))$ in (7.14) and $N = \mathbf{J}(1+|R|)$, one gets (7.3). □

Let $r_{\min} > 0$ be the minimum $r$ such that (7.1)–(7.3) are satisfied for some $\epsilon > 0$. Now for each $r > r_{\min}$, we let $\epsilon(r)$ be the smallest $\epsilon$ such that (7.1)–(7.3) are satisfied. By Proposition 1, it is easy to see that

$$\lim_{r\to\infty} \epsilon(r) = 0.$$

With $\epsilon(r)$ replacing $\epsilon$, we call each of the events in (7.1)–(7.3) a "bad" event. We let $\mathcal{G}^r$ denote the intersection of the complements of these "bad" events, and the events in $\mathcal{G}^r$ are referred to as "good" events in the rest of this section. Obviously,

$$\lim_{r\to\infty} \mathbb{P}(\mathcal{G}^r) = 1.$$

7.2. *State space collapse on* $[\xi_{r,0}\tau_0/r^2, T]$. Again, let $L > 1 \vee \tau_0$ and $T > 0$ be fixed. We divide the diffusion-scaled time interval $[0, T]$ into two overlapping intervals: $[0, \xi_{r,0}L/r^2]$ and $[\xi_{r,0}\tau_0/r^2, T]$. In this section we show a state space collapse result on the time interval $[\xi_{r,0}\tau_0/r^2, T]$. The state space collapse on $[0, \xi_{r,0}L/r^2]$ will be presented in the next section.

In order to connect the fluid-scaled processes with the fluid model, we first introduce the notion of *cluster point*. Let $F = \mathbb{D}^d[0, L]$ be the space of right continuous functions with left limits from $[0, L]$ to $\mathbb{R}^{\mathbf{I}+\mathbf{J}}$, where $d = \mathbf{I} + \mathbf{J}$. Let $\mathcal{C} = \{F_r\}$ be a sequence of subsets of $F$. A point $f$ in $F$ is said to be a *cluster point* of $\mathcal{C}$ if for each $\epsilon > 0$ and $r_0 > 0$, there exist $r \geq r_0$ and $g \in F_r$ such that $f - g_L < \epsilon$. The sequence $\{F_r\}$ is said to be *asymptotically Lipschitz* if there exist a constant $\kappa > 0$ and a sequence of positive numbers $\{\epsilon(r)\}$ with $\epsilon(r) \to 0$ as $r \to \infty$ such that for each $r$ all elements $f \in F_r$ satisfy both

$$|f(0)| \leq 1 \tag{7.15}$$

and

$$|f(t_2) - f(t_1)| \leq \kappa |t_2 - t_1| + \epsilon(r) \qquad \text{for all } t_1, t_2 \in [0, L]. \tag{7.16}$$

Let $F'$ denote those $f \in F$ satisfying both (7.15) and

$$|f(t_2) - f(t_1)| \leq \kappa |t_2 - t_1| \qquad \text{for all } t_1, t_2 \in [0, L]. \tag{7.17}$$

The following lemma is due to Bramson (1998). We state it here for completeness.

LEMMA 7 [Bramson (1998), Proposition 4.1]. *Assume that $\mathcal{C}$ is asymptotically Lipschitz. For each $\epsilon > 0$, there exists an $r_0$, so that for each $r \geq r_0$ and $g \in F_r$, one has $\|f - g\|_L < \epsilon$ for some cluster point $f$ of $\mathcal{C}$ with $f \in F'$.*

Lemma 7 says that the cluster points are "rich": for large $r$, all elements in $F_r$ can be uniformly approximated by cluster points.

We set

$$F_g^r = \{\mathbb{X}^{r,m}(\cdot, \omega), m < rT, \omega \in \mathcal{G}^r\} \qquad \text{for each } r$$

and

$$\mathcal{F}_g = \{F_g^r\}. \tag{7.18}$$

From the choice of our fluid scale $\xi_{r,m}$ and the definition of $T^{r,m}(\cdot)$, it follows that

$$|\mathbb{X}^{r,m}(0)| \le 1.$$

It follows from (7.3) in Proposition 1 that the sequence of sets of scaled stochastic processing network processes $\mathbb{X}^{r,m}(\cdot)$ is asymptotically Lipschitz. Lemma 7 immediately implies the following proposition which says that, for large $r$, $\mathbb{X}^{r,m}(\cdot)$ are uniformly close to cluster points that are Lipschitz continuous.

PROPOSITION 2. *Fix $\epsilon > 0$, $L$ and $T$, and choose $r$ large enough. Then, for $\omega \in \mathcal{G}^r$ and any $m < rT$,*

$$\|\mathbb{X}^{r,m}(\cdot,\omega) - \tilde{\mathbb{X}}(\cdot)\|_L \le \epsilon$$

*for some cluster point $\tilde{\mathbb{X}}(\cdot)$ of $\mathcal{F}_g$ with $\tilde{\mathbb{X}}(\cdot) \in F'$.*

The next proposition says that if the stochastic processing networks operate under a maximum pressure policy, then each cluster point of $\mathcal{F}_g$ satisfies fluid model equations (6.1)–(6.6), and thus is a fluid model solution to the fluid model operating under the maximum pressure policy.

PROPOSITION 3. *Consider a sequence of stochastic processing networks operating under a maximum pressure policy. Assume Assumptions 1 and 3 and that the limit network satisfies the EAA condition. Fix $L > 0$ and $T > 0$. Then all cluster points of $\mathcal{F}_g$ are solutions to the fluid model equations (6.1)–(6.6) on $[0, L]$, and they are all in $F'$.*

PROOF. The idea is to approximate each cluster point of $\mathcal{F}_g$ with some $\mathbb{X}^{r,m}$ on $\mathcal{G}^r$ and show that the equations (6.1)–(6.6) are asymptotically satisfied by $\mathbb{X}^{r,m}$. We will only demonstrate equations (6.1) and (6.6); the others are quite straightforward and can be verified similarly. Let $\tilde{\mathbb{X}} = (\tilde{Z}, \tilde{T})$ be a cluster point. We first verify equation (6.1). For any $\epsilon > 0$ and all $r$ that are large enough, there are $\omega \in \mathcal{G}^r$ and $m < rT$ such that

$$\|S_j^{r,m}(T_j^{r,m}(t)) - \mu_j^r T_j^{r,m}(t)\|_L \le \epsilon \qquad \text{for each } j \in \mathcal{J}, \tag{7.19}$$

$$\|\Phi_i^{j,r,m}(S_j^{r,m}(T_j^{r,m}(t))) - P_i^j \mu_j^r T_j^{r,m}(t)\|_L \le \epsilon \qquad \text{for each } j \in \mathcal{J}, i \in \mathcal{B}_j, \tag{7.20}$$

$$\|\tilde{Z}(\cdot) - Z^{r,m}(\cdot)\|_L \le \epsilon, \tag{7.21}$$

$$\|\tilde{T}(\cdot) - T^{r,m}(\cdot)\|_L \le \epsilon, \tag{7.22}$$

$$|R^r - R| \le \epsilon. \tag{7.23}$$

Scaling (5.2) and plugging in the bounds in (7.19) and (7.21), we have

$$|Z^{r,m}(t) - Z^{r,m}(0) + R^r T^{r,m}(t)| \le 2\epsilon \qquad \text{for all } t \in [0, L].$$

From (7.22) and (7.23), we have, for each $t \in [0, L]$,

$$\begin{aligned} |R^r T^{r,m}(t) - R\tilde{T}(t)| &\le |R^r(T^{r,m}(t) - \tilde{T}(t))| + |(R^r - R)\tilde{T}(t)| \\ &\le N\mathbf{J}\epsilon + \epsilon N L \mathbf{J}. \end{aligned}$$

Recall that $N \ge \sup_r |R^r|$. It then follows that, for each $t \le L$,

$$\begin{aligned} |\tilde{Z}(t) - \tilde{Z}(0) + R\tilde{T}(t)| \le\ & |\tilde{Z}(t) - Z^{r,m}(t)| + |Z^{r,m}(0) - \tilde{Z}(0)| \\ & + |R^r T^{r,m}(t) - R\tilde{T}(t)| \\ & + |Z^{r,m}(t) - Z^{r,m}(0) + R^r T^{r,m}(t)| \\ \le\ & (4 + NL\mathbf{J} + N\mathbf{J})\epsilon. \end{aligned}$$

Thus, equation (6.1) is satisfied by $\tilde{\mathbb{X}}$ because $\epsilon$ can be arbitrarily small.

To show equation (6.6), first observe that, for any allocation $a \in \mathcal{A}$,

$$\begin{aligned} |p_\alpha(a, \tilde{Z}(t)) - p^r_\alpha(a, Z^{r,m}(t))| &= |\alpha \times \tilde{Z}(t) \cdot Ra - \alpha \times Z^{r,m}(t) \cdot R^r a| \\ &\le |\alpha|(|\tilde{Z}(t) \cdot (Ra - R^r a)| \\ &\qquad + |(\tilde{Z}(t) - Z^{r,m}(t)) \cdot R^r a|) \\ &\le |\alpha|((NL+1)\mathbf{IJ}\epsilon + \epsilon \mathbf{IJ} N) \\ &= (NL + N + 1)|\alpha|\mathbf{IJ}\epsilon. \end{aligned} \tag{7.24}$$

In (7.24), $p^r_\alpha(a, q) \equiv (\alpha \times q) \cdot R^r a$ is the network pressure for the $r$th network, associated with the parameter $\alpha$, under allocation $a$ when the queue length is $q$.

Denote $\mathcal{E}^* = \arg\max_{a \in \mathcal{E}} p_\alpha(a, \tilde{Z}(t))$ as the set of maximum extreme allocations under buffer size $\tilde{Z}(t)$. Because the limit network satisfies the EAA condition, we can choose an $a^* \in \mathcal{E}^*$ such that $\tilde{Z}_i(t) > 0$ for each constituent buffer $i$ of $a^*$. Denote $\mathcal{I}(a^*)$ the set of constituent buffers. Namely,

$$\mathcal{I}(a^*) = \left\{ i : \sum_j a^*_j B_{ji} > 0 \right\}.$$

Then $\tilde{Z}_i(t) > 0$ for all $i \in \mathcal{I}(a^*)$.

Suppose $a \in \mathcal{E} \setminus \mathcal{E}^*$. Then $p_\alpha(a, \tilde{Z}(t)) < \max_{a' \in \mathcal{E}} p_\alpha(a', \tilde{Z}(t))$. Since

$$p_\alpha(a^*, \tilde{Z}(t)) = \max_{a' \in \mathcal{E}} p_\alpha(a', \tilde{Z}(t))$$

and $\min_{i \in \mathcal{I}(a^*)} \tilde{Z}_i(t) > 0$, by the continuity of $\tilde{\mathbb{X}}(\cdot)$, there exist $\epsilon_1 > 0$ and $\delta > 0$ such that, for each $\tau \in [t - \delta, t + \delta]$ and $i \in \mathcal{I}(a^*)$,

$$p_\alpha(a, \tilde{Z}(\tau)) + \epsilon_1 \le p_\alpha(a^*, \tilde{Z}(\tau)) \quad \text{and} \quad \tilde{Z}_i(\tau) \ge \epsilon_1. \tag{7.25}$$

For sufficiently large $r$, we can choose $\epsilon$ small enough such that $(NL + N + 1) \times |\alpha|\mathbf{IJ}\epsilon \le \epsilon_1/3$. It follows from (7.21), (7.24) and (7.25) that, for all $\tau \in [t - \delta, t + \delta]$,

$$p^r_\alpha(a, Z^{r,m}(\tau)) + \epsilon_1/3 \le p^r_\alpha(a^*, Z^{r,m}(\tau)),$$
$$Z^{r,m}_i(\tau) \ge \epsilon_1/2 \qquad \text{for each } i \in \mathcal{I}(a^*).$$

Choosing $r > 2\mathbf{J}/\epsilon_1$, for each $\tau \in [rm + \xi_{r,m}(t-\delta), rm + \xi_{r,m}(t+\delta)]$, we have

$$p^r_\alpha(a, Z^r(\tau)) < p^r_\alpha(a^*, Z^r(\tau)), \tag{7.26}$$
$$Z^r_i(\tau) \ge \mathbf{J} \qquad \text{for each } i \in \mathcal{I}(a^*). \tag{7.27}$$

Condition (7.27) implies that $a^*$ is a feasible allocation at any time $\tau \in [rm + \xi_{r,m}(t-\delta), rm + \xi_{r,m}(t+\delta)]$, that is, $a^* \in \mathcal{E}(\tau)$. Following (7.26) and the definition of a (preemptive-resume) maximum pressure policy, the allocation $a$ will not be employed during time interval $[rm + \xi_{r,m}(t-\delta), rm + \xi_{r,m}(t+\delta)]$. Therefore, only the allocations in $\mathcal{E}^*$ will be employed during this interval.

For each $a \in \mathcal{E}$, denote $(T^a)^r(t)$ to be the cumulative amount of time allocation $a$ has been employed by time $t$. Because only allocations in $\mathcal{E}$ are employed under a maximum pressure policy, we have that, for each $r$ and all $t \ge 0$, under a maximum pressure policy,

$$\sum_{a\in\mathcal{E}} (T^a)^r(t) = t.$$

Furthermore, since the employment of allocation $a$ for one unit of time contributes to $a_j$ unit of activity $j$ processing time, we have

$$T^r(t) = \sum_{a\in\mathcal{E}} a(T^a)^r(t).$$

Then it follows that

$$\begin{aligned}
&(\alpha \times \tilde{Z}(t)) \cdot (R[T^r(rm + \xi_{r,m}(t+\delta)) - T^r(rm + \xi_{r,m}(t-\delta))]) \\
&\quad = (\alpha \times \tilde{Z}(t)) \cdot \Bigg(R\Bigg[\sum_{a\in\mathcal{E}} a((T^a)^r(rm + \xi_{r,m}(t+\delta)) \\
&\qquad\qquad - (T^a)^r(rm + \xi_{r,m}(t-\delta)))\Bigg]\Bigg) \\
&\quad = \sum_{a\in\mathcal{E}} ((\alpha \times \tilde{Z}(t)) \cdot Ra)[(T^a)^r(rm + \xi_{r,m}(t+\delta)) \\
&\qquad\qquad - (T^a)^r(rm + \xi_{r,m}(t-\delta))] \\
&\quad = \sum_{a\in\mathcal{E}^*} ((\alpha \times \tilde{Z}(t)) \cdot Ra)[(T^a)^r(rm + \xi_{r,m}(t+\delta))
\end{aligned}$$

(7.28)

$$
\begin{aligned}
&\qquad\qquad - (T^a)^r(rm + \xi_{r,m}(t-\delta))] \\
&= \Big(\max_{a\in\mathcal{E}}(\alpha \times \tilde{Z}(t)) \cdot Ra\Big) \sum_{a\in\mathcal{E}^*} [(T^a)^r(rm + \xi_{r,m}(t+\delta)) \\
&\qquad\qquad - (T^a)^r(rm + \xi_{r,m}(t-\delta))] \\
&= \Big(\max_{a\in\mathcal{E}}(\alpha \times \tilde{Z}(t)) \cdot Ra\Big) \sum_{a\in\mathcal{E}} [(T^a)^r(rm + \xi_{r,m}(t+\delta)) \\
&\qquad\qquad - (T^a)^r(rm + \xi_{r,m}(t-\delta))] \\
&= 2\xi_{r,m}\delta\Big(\max_{a\in\mathcal{E}}(\alpha \times \tilde{Z}(t)) \cdot Ra\Big).
\end{aligned}
$$

The second and fourth equalities in (7.28) hold because only allocations in $\mathcal{E}^*$ will be employed during $[rm + \xi_{r,m}(t-\delta), rm + \xi_{r,m}(t+\delta)]$; the third holds because every allocation $a \in \mathcal{E}^*$ has the same network pressure equal to $\max_{a\in\mathcal{E}}(\alpha \times \tilde{Z}(t)) \cdot Ra$. From (7.28), we have

$$(\alpha \times \tilde{Z}(t)) \cdot R(T^{r,m}(t+\delta) - T^{r,m}(t-\delta))/2\delta = \max_{a\in\mathcal{E}}(\alpha \times \tilde{Z}(t)) \cdot Ra.$$

Because $\epsilon$ in (7.22) can be arbitrarily small, we have

$$(\alpha \times \tilde{Z}(t)) \cdot R(\tilde{T}(t+\delta) - \tilde{T}(t-\delta))/2\delta = \max_{a\in\mathcal{E}}(\alpha \times \tilde{Z}(t)) \cdot Ra,$$

and by letting $\delta \to \infty$, (6.6) is verified.

Since $\tilde{\mathbb{X}}$ and $\mathbb{X}^{r,m}$ are uniformly close and $\mathbb{X}^{r,m}$ satisfies (7.15) and (7.16), it is straightforward to verify that $\tilde{\mathbb{X}}$ satisfies both (7.15) and (7.17), thus showing $\tilde{\mathbb{X}} \in F'$. $\square$

Since every fluid model solution under the maximum pressure policy satisfies (6.7) by Theorem 6, Proposition 3 implies that any cluster point $\tilde{\mathbb{X}}$ of $\mathcal{F}_g$ satisfies (6.7). That is,

$$|\tilde{Z}(t) - \zeta^\alpha \tilde{W}(t)| = 0 \qquad \text{for } t \geq \tau_0,$$

where $\tilde{W}(t) = y^* \cdot \tilde{Z}(t)$. Define

$$W^{r,m} = y^r \cdot Z^{r,m}.$$

Because the fluid-scaled stochastic processing network processes can be uniformly approximated by cluster points, it leads to the following proposition.

PROPOSITION 4. *Fix $L > \tau_0$, $T > 0$ and $\epsilon > 0$. For $r$ large enough,*

$$|Z^{r,m}(t) - \zeta^\alpha W^{r,m}(t)| \leq \epsilon \qquad \textit{for all } 0 \leq m \leq rT, \tau_0 \leq t \leq L, \omega \in \mathcal{G}^r.$$

PROOF. From Proposition 2, for each $\epsilon > 0$, large enough $r$, and each $0 \le m \le rT$ and $\omega \in \mathcal{G}^r$, we can find a cluster point $\tilde{\mathbb{X}}$ such that

$$\|\tilde{Z}(\cdot) - Z^{r,m}(\cdot)\|_L \le \epsilon.$$

From Proposition 3 and Theorem 6, we have

$$|\tilde{Z}(t) - \zeta^\alpha \tilde{W}(t)| = 0 \qquad \text{for } t \ge \tau_0.$$

Then, for each $t \ge \tau_0$,

$$|Z^{r,m}(t) - \zeta^\alpha W^{r,m}(t)| \le |Z^{r,m}(t) - \tilde{Z}(t)| + |\zeta^\alpha \tilde{W}(t) - \zeta^\alpha W^{r,m}(t)|. \tag{7.29}$$

Because $y^r \to y^*$ as $r \to \infty$, for large enough $r$,

$$|y^r - y^*| < \epsilon.$$

Let $\kappa = (\sup_r |y^r|) \vee |\zeta^\alpha| < \infty$. Then we have

$$\begin{aligned} |\tilde{W}(t) - W^{r,m}(t)| &\le |y^* \cdot \tilde{Z}(t) - y^r \cdot \tilde{Z}(t)| + |y^r \cdot \tilde{Z}(t) - y^r \cdot Z^{r,m}(t)| \\ &\le \mathbf{I}(\epsilon|\tilde{Z}(t)| + \kappa\epsilon) \le (NL + 1 + \kappa)\mathbf{I}\epsilon. \end{aligned}$$

The last inequality follows from the fact that $\tilde{Z}(t) \le NL + 1$ for all $t \le L$. Note that $\tilde{Z}(t) \in F'$ and thus satisfies (7.17). One also gets $\tilde{W}(t) \le \kappa(NL+1)$ since $|y^*| \le \kappa$. From (7.29), we have

$$|Z^{r,m}(t) - \zeta^\alpha W^{r,m}(t)| \le \epsilon + \kappa(NL + 1 + \kappa)\mathbf{I}\epsilon.$$

The proposition follows by rechoosing $\epsilon$. $\square$

We need to translate the results in Proposition 4 into the state space collapse results for diffusion-scaled processes. First we can express $Z^{r,m}(t)$ by $\hat{Z}^r$ via

$$Z^{r,m}(t) = \frac{r}{\xi_{r,m}} \hat{Z}^r\left(\frac{t\xi_{r,m} + rm}{r^2}\right) = \frac{1}{\bar{\xi}_{r,m}} \hat{Z}^r\left(\frac{t\bar{\xi}_{r,m} + m}{r}\right),$$

where $\bar{\xi}_{r,m} = \xi_{r,m}/r$. The time interval $[\tau_0, L]$ for the fluid-scaled process $Z^{r,m}$ corresponds to the time interval $[(m + \bar{\xi}_{r,m}\tau_0)/r, (m + \bar{\xi}_{r,m}L)/r]$ for the diffusion-scaled process $\hat{Z}^r$. Proposition 4 immediately leads to the following.

PROPOSITION 5. *Fix $L > \tau_0$, $T > 0$ and $\epsilon > 0$. For $r$ large enough and each $m < rT$,*

$$\begin{aligned} &|\widehat{Z}^r(t) - \zeta^\alpha \widehat{W}^r(t)| \le \bar{\xi}_{r,m}\epsilon \\ &\textit{for all } (m + \bar{\xi}_{r,m}\tau_0)/r \le t \le (m + \bar{\xi}_{r,m}L)/r,\ \omega \in \mathcal{G}^r. \end{aligned} \tag{7.30}$$

Proposition 5 gives estimates on each small interval for $|\widehat{Z}^r(t) - \zeta^\alpha \widehat{W}^r(t)|$. We shall obtain the estimate on the whole time interval $[0,T]$, and then show that $\bar{\xi}_{r,m}$ are stochastically bounded. The following lemma ensures that for large enough $L$, in particular, for $L \geq 3N\tau_0 + 1$, the small intervals in Proposition 5 are overlapping, and therefore, the estimate on $[\bar{\xi}_{r,0}\tau_0/r, T]$ is obtained.

LEMMA 8. *For a fixed $T > 0$ and large enough $r$,*

$$\bar{\xi}_{r,m+1} \leq 3N\bar{\xi}_{r,m}$$

*for $\omega \in \mathcal{G}^r$ and $m < rT$, where $N$ is chosen as in Proposition 1.*

PROOF. By the definition of $\mathcal{G}^r$,

$$|Z^{r,m}(t_2) - Z^{r,m}(t_1)| \leq N|t_2 - t_1| + 1$$

for $t_1, t_2 \in [0, L]$ and $m < rT$. Setting $t_1 = 0$ and $t_2 = 1/\bar{\xi}_{r,m}$, we have

$$|\hat{Z}^r((m+1)/r)| - |\hat{Z}^r(m/r)| \leq N + \bar{\xi}_{r,m}. \tag{7.31}$$

From the definition of $\bar{\xi}_{r,m+1}$, we have

$$\begin{aligned}\bar{\xi}_{r,m+1} &= |\hat{Z}^r((m+1)/r)| \vee 1 \\ &\leq (|\hat{Z}^r(m/r)| + N + \bar{\xi}_{r,m}) \vee 1 \\ &\leq N + 2\bar{\xi}_{r,m} \leq 3N\bar{\xi}_{r,m}.\end{aligned}$$

The first inequality follows from (7.31), and the second inequality follows from the definition of $N$ and $\bar{\xi}_{r,m}$. □

7.3. *State space collapse on* $[0, \xi_{r,0}L/r^2]$. Now we shall estimate $|\widehat{Z}^r(t) - \zeta^\alpha \widehat{W}^r(t)|$ on the interval $[0, \bar{\xi}_{r,0}L/r] = [0, \xi_{r,0}L/r^2]$. This will be given by the initial condition (4.10) and the result in the second part of Theorem 6.

Condition (4.10) implies that

$$Z^{r,0}(0) \to 0 \qquad \text{in probability.}$$

Then, for each $r > 0$, we let

$$\epsilon_1(r) = \min\{\epsilon : \mathbb{P}(Z^{r,0}(0) > \epsilon) \leq \epsilon\}. \tag{7.32}$$

It follows from (7.32) that

$$\lim_{r\to\infty} \epsilon_1(r) \to 0.$$

Now let $\mathcal{L}^r$ be the intersection of $\mathcal{G}^r$ and the event

$$Z^{r,0}(0) \leq \epsilon_1(r).$$

Obviously, $\lim_{r\to\infty} \mathbb{P}(\mathcal{L}^r) = 1$.

We define

$$\mathcal{F}_o = \{F_o^r\}$$

with

$$F_o^r = \{\mathbb{X}^{r,0}(\cdot,\omega), \omega \in \mathcal{L}^r\}.$$

Parallel to Proposition 2, we have the following proposition which states that $\mathcal{F}_o$ can be asymptotically approximated by cluster points of $\mathcal{F}_o$.

PROPOSITION 6. *Fix $\epsilon > 0$, $L > 0$ and $T > 0$, and choose $r$ large enough. Then, for $\omega \in \mathcal{L}^r$,*

$$\|\mathbb{X}^{r,0}(\cdot,\omega) - \tilde{\mathbb{X}}(\cdot)\|_L \le \epsilon$$

*for some cluster point $\tilde{\mathbb{X}}(\cdot)$ of $\mathcal{F}_o$ with $\tilde{\mathbb{X}}(\cdot) \in F'$.*

PROOF. Since both (7.15) and (7.16) are satisfied by $\mathbb{X}^{r,0}$, the result follows from Lemma 7. □

PROPOSITION 7. *Fix $L > 0$. Then for any cluster point $\tilde{\mathbb{X}}(\cdot)$ of $\mathcal{F}_o$,*

$$\tilde{Z}(t) = \zeta^\alpha \tilde{W}(t) \qquad \text{for } t \in [0, L].$$

PROOF. Since any cluster point of $\mathcal{F}_o$ is automatically a cluster point of $\mathcal{F}_g$, therefore, it satisfies all the fluid model equations. It suffices to show that

$$|\tilde{Z}(0) - \zeta^\alpha \tilde{W}(0)| = 0, \tag{7.33}$$

which, together with Theorem 6, implies the result. In fact, we have $\tilde{Z}(0) = 0$, since for any given $\delta > 0$, one can choose $r$ large enough and $\omega \in \mathcal{L}^r$ such that $|\tilde{Z}(0) - Z^{r,0}(0)| \le \delta$, and $\mathbb{Z}^{r,0}(0) \le \delta$. □

Propositions 6 and 7 immediately lead to the following proposition, which is parallel to Propositions 4 and 5.

PROPOSITION 8. *Fix $L$ and $\epsilon > 0$. For large enough $r$,*

$$|Z^{r,0}(t) - \zeta^\alpha W^{r,0}(t)| \le \epsilon \qquad \text{for all } 0 \le t \le L, \omega \in \mathcal{L}^r$$

*and*

$$|\widehat{Z}^r(t) - \zeta^\alpha \widehat{W}^r(t)| \le \bar{\xi}_{r,0}\epsilon \qquad \text{for all } 0 \le t \le \bar{\xi}_{r,0}L/r, \omega \in \mathcal{L}^r. \tag{7.34}$$

7.4. *Proof of Theorem* 5. Propositions 5 and 8 together give the multiplicative state space collapse of the stochastic processing network processes. To prove the state space collapse result stated in Theorem 5, it is enough to prove that $\bar{\xi}_{r,m}$ are stochastically bounded. We first give an upper bound on $\bar{\xi}_{r,m}$ in terms of $\widehat{W}^r$.

LEMMA 9. *If* $|\widehat{Z}^r(m/r) - \zeta^\alpha \widehat{W}^r(m/r)| \le 1$, *there exists some* $\kappa \ge 1$ *such that*

$$\bar{\xi}_{r,m} \le \kappa(\widehat{W}^r(m/r) \vee 1).$$

PROOF. Because $\xi_{r,m} = |Z^r(rm)| \vee r$, we have

$$\bar{\xi}_{r,m} = |\widehat{Z}^r(m/r)| \vee 1 \le |\zeta^\alpha|\widehat{W}^r(m/r) + 1 \le 2(|\zeta^\alpha|\widehat{W}^r(m/r) \vee 1).$$

The result then follows by choosing $\kappa = 2(|\zeta^\alpha| \vee 1)$. □

The following proposition will be used to derive an upper bound on the oscillation of $\widehat{W}^r$.

PROPOSITION 9. *Consider a sequence of stochastic processing networks operating under the maximum pressure policy with parameter* $\alpha$. *There exists* $\epsilon_0 > 0$ *such that, for large enough* $r$, *and any* $0 \le t_1 < t_2$, *if* $\widehat{W}^r(t) \ge \mathbf{J}/(r\epsilon_0)$ *and* $|\widehat{Z}^r(t)/\widehat{W}^r(t) - \zeta^\alpha| \le \epsilon_0$ *for all* $t \in [t_1, t_2]$, *then*

$$\widehat{Y}^r(t_2) = \widehat{Y}^r(t_1).$$

PROOF. First, we can choose $0 < \epsilon_0 \le \min\{\zeta_i^\alpha/2 : \zeta_i^\alpha > 0\}$ so that for all $i$ with $\zeta_i^\alpha > 0$ and all $t \in [t_1, t_2]$, we have

$$\widehat{Z}_i^r(t) \ge \widehat{W}^r(t)(\zeta_i^\alpha - \epsilon_0) \ge \widehat{W}^r(t)\epsilon_0 \ge \mathbf{J}/r.$$

That is, for all $i$ with $\zeta_i^\alpha > 0$,

$$Z_i^r(\tau) \ge \mathbf{J} \qquad \text{for all } \tau \in [r^2 t_1, r^2 t_2]. \tag{7.35}$$

Define

$$\mathcal{E}^* = \arg\max_{a \in \mathcal{E}} y^* \cdot Ra = \arg\max_{a \in \mathcal{E}} (\alpha \times \zeta^\alpha) \cdot Ra.$$

The second equality follows from the fact that $\alpha \times \zeta^\alpha = y^* / \sum_{i \in \mathcal{I}}((y_i^*)^2/\alpha_i)$.

We now show that at least one of the allocations in $\mathcal{E}^*$ is feasible during $[r^2 t_1, r^2 t_2]$. Define $\mathcal{I}(a) = \{i \in \mathcal{I} : \sum_j B_{ji} a_j > 0\}$ to be the constituency buffers of the allocation $a$. One sufficient condition for an allocation $a$ to be feasible at time $t$ is that $Z_i^r(t) \ge \mathbf{J}$ for all $i \in \mathcal{I}(a)$. Because the limit network satisfies the EAA condition, replacing $q$ by $y^*$ in Definition 4 implies that there exists an allocation $a^* \in \mathcal{E}^*$ such that $y_i^* > 0$ for each $i \in \mathcal{I}(a^*)$; namely, all constituency buffers of $a^*$ have positive $y_i^*$'s. Since $y_i^* > 0$ implies $\zeta_i^\alpha > 0$, it follows from (7.35) that $Z_i^r(\tau) \geq \mathbf{J}$ for all $\tau \in [r^2 t_1, r^2 t_2]$ and $i \in \mathcal{I}(a^*)$. This implies that $a^*$ is a feasible allocation during $[r^2 t_1, r^2 t_2]$.

Next, we will show that if $\epsilon_0$ is chosen sufficiently small, only allocations in $\mathcal{E}^*$ can be employed under a maximum pressure policy during $[r^2 t_1, r^2 t_2]$. Let

$$\epsilon_1 = \left( y^* \cdot Ra^* - \max_{a \in \mathcal{E} \setminus \mathcal{E}^*} y^* \cdot Ra \right) > 0.$$

Set

$$\kappa_0 = \sup_r \max_{a \in \mathcal{E}} |R^r a| \quad \text{and} \quad \kappa_1 = \sum_{i \in \mathcal{I}} ((y_i^*)^2 / \alpha_i).$$

Choose $r$ large enough such that $|(y^*)'(R^r - R)| \leq \epsilon_0$. Then for each allocation $a \in \mathcal{E} \setminus \mathcal{E}^*$ and each $\tau \in [r^2 t_1, r^2 t_2]$,

$$\begin{aligned}
&((\alpha \times Z^r(\tau)) \cdot R^r a^* - (\alpha \times Z^r(\tau)) \cdot R^r a)/W^r(\tau) \\
&\quad = \left( \frac{\alpha \times Z^r(\tau)}{W^r(\tau)} - \alpha \times \zeta^\alpha \right) \cdot R^r a^* \\
&\qquad + (y^* \cdot R^r a^* - y^* \cdot Ra^*)/\kappa_1 + (y^* \cdot Ra^* - y^* \cdot Ra)/\kappa_1 \\
&\qquad + (y^* \cdot Ra - y^* \cdot R^r a)/\kappa_1 + \left( \alpha \times \zeta^\alpha - \frac{\alpha \times Z^r(\tau)}{W^r(\tau)} \right) \cdot R^r a \\
&\quad \geq -|\alpha| \epsilon_0 \kappa_0 \mathbf{I} - \epsilon_0 \mathbf{J}/\kappa_1 + \epsilon_1/\kappa_1 - \epsilon_0 \mathbf{J}/\kappa_1 - |\alpha| \epsilon_0 \kappa_0 \mathbf{I} \\
&\quad = (\epsilon_1 - 2(|\alpha| \kappa_0 \kappa_1 \mathbf{I} + \mathbf{J}) \epsilon_0)/\kappa_1.
\end{aligned} \tag{7.36}$$

Thus,

$$p_\alpha^r(a^*, Z^r(\tau)) - p_\alpha^r(a, Z^r(\tau)) \geq W^r(\tau)(\epsilon_1 - 2(|\alpha| \kappa_0 \kappa_1 \mathbf{I} + \mathbf{J})\epsilon_0)/\kappa_1.$$

Set $\epsilon_0 = \epsilon_1/(3(|\alpha|)\kappa_0 \kappa_1 \mathbf{I} + \mathbf{J})) \wedge \min_{i \in \mathcal{I}} \{\zeta_i^\alpha / 2 : \zeta_i^\alpha > 0\}$, then

$$p_\alpha^r(a^*, Z^r(\tau)) - p_\alpha^r(a, Z^r(\tau)) > 0. \tag{7.37}$$

Obviously, $\epsilon_0$ depends only on $\alpha, R$ and $A$. Equation (7.37) implies that the pressure under allocation $a^*$ is strictly larger than that under allocation $a$. It follows that only the allocations in $\mathcal{E}^*$ can be employed during $[r^2 t_1, r^2 t_2]$. Therefore, for large enough $r$,

$$\begin{aligned}
y^r \cdot R^r (T^r(r^2 t_2) - T^r(r^2 t_1)) &= \sum_{a \in \mathcal{E}} y^r \cdot R^r a((T^a)^r(r^2 t_2) - (T^a)^r(r^2 t_1)) \\
&= \sum_{a \in \mathcal{E}^*} y^r \cdot R^r a((T^a)^r(r^2 t_2) - (T^a)^r(r^2 t_1))
\end{aligned} \tag{7.38}$$

$$= \sum_{a \in \mathcal{E}^*} (1-\rho^r)((T^a)^r(r^2 t_2) - (T^a)^r(r^2 t_1))$$
$$= (1-\rho^r)(r^2 t_2 - r^2 t_1),$$

where the third equality follows from Lemma 10 below. Equations (7.38) and (5.1) imply that $\widehat{Y}^r(t_2) - \widehat{Y}^r(t_1) = 0$. □

LEMMA 10. *For each $a^* \in \mathcal{E}^* = \arg\max_{a \in \mathcal{E}} y^* \cdot Ra$ and large enough $r$,*

$$y^r \cdot R^r a^* = \max_{a \in \mathcal{E}} y^r \cdot R^r a = 1 - \rho^r.$$

The proof of Lemma 10 will be provided in Appendix B. We now complete the proof of Theorem 5.

PROOF OF THEOREM 5. We only need to show that $\bar{\xi}_{r,m}$ are stochastically bounded. We first fix an $\epsilon > 0$. Because $\hat{W}^r(0)$ is stochastically bounded, there exists a $\kappa_1$ such that, for $r$ large enough,

$$\mathbb{P}(\hat{W}^r(0) \le \kappa_1) \ge 1 - \epsilon.$$

Recall that the process $X^*$, defined in Lemma 3, is a Brownian motion, so it has continuous sample path almost surely. Therefore, there exists a $\kappa_2$ such that

$$\mathbb{P}(\mathrm{Osc}(X^*, [0, T+L]) \le \kappa_2/2) \ge 1 - \epsilon/2,$$

where

$$\mathrm{Osc}(f, [0,t]) = \sup_{0 \le t_1 < t_2 \le t} f|(t_2) - f(t_1)|.$$

Since $\widehat{X}^r$ converges to $X^*$ in distribution, for large enough $r$,

$$\mathbb{P}(\mathrm{Osc}(\hat{X}^r, [0, T+L]) \le \kappa_2) \ge 1 - \epsilon.$$

Meanwhile, because

$$|\widehat{Z}^r(0) - \zeta^\alpha \widehat{W}^r(0)| \to 0 \qquad \text{in probability},$$

we have

$$\mathbb{P}(|\widehat{Z}^r(0) - \zeta^\alpha \widehat{W}^r(0)| \le \epsilon) \ge 1 - \epsilon$$

for $r$ large enough.

Define

$$\mathcal{H}^{r,\epsilon} = \{\omega : \widehat{W}^r(0) \le \kappa_1, \mathrm{Osc}(\hat{X}^r, [0, T+L]) \le \kappa_2,$$
$$\text{and } |\widehat{Z}^r(0) - \zeta^\alpha \widehat{W}^r(0)| \le \epsilon\}.$$

Then for $r$ large enough,

$$\mathbb{P}(\mathcal{H}^{r,\epsilon}) \geq 1 - 3\epsilon.$$

Furthermore, we can choose $r$ large enough such that Propositions 5 and 8 hold with $\epsilon$ replaced by $\epsilon/\kappa(\kappa_1+\kappa_2+1)$ and $\mathbb{P}(\mathcal{L}^r) \geq 1-\epsilon$, where $\kappa$ is given as in Lemma 9. Note that $\mathbb{P}(\mathcal{L}^r) \to 1$ as $r \to \infty$.

Denote $\mathcal{N}^{r,\epsilon} = \mathcal{L}^r \cap \mathcal{H}^{r,\epsilon}$, then for all $r$ large enough,

$$\mathbb{P}(\mathcal{N}^{r,\epsilon}) \geq 1 - 4\epsilon.$$

Now we are going to show that if $\epsilon \leq \epsilon_0$, $\bar{\xi}_{r,m} \leq \kappa(\kappa_1+\kappa_2+1)$ on $\mathcal{N}^{r,\epsilon}$ for all $r$ large enough and $m \leq rT$. In fact, we are going to show the following is true on $\mathcal{N}^{r,\epsilon}$ for all $r$ large enough and $m \leq rT$:

$$|\widehat{Z}^r(m/r) - \zeta^\alpha \widehat{W}^r(m/r)| \leq \epsilon, \tag{7.39}$$

$$\bar{\xi}_{r,m} \leq \kappa(\kappa_1+\kappa_2+1), \tag{7.40}$$

$$|\widehat{Z}^r(t) - \zeta^\alpha \widehat{W}^r(t)| \leq \epsilon \qquad \text{for all } 0 \leq t \leq (m+\bar{\xi}_{r,m}L)/r, \tag{7.41}$$

$$\int_0^{(m+\bar{\xi}_{r,m}L)/r} 1_{(\widehat{W}^r(s)>1)}\, d\widehat{Y}^r(s) = 0. \tag{7.42}$$

This will be shown by induction. When $m=0$, (7.39) obviously holds on $\mathcal{N}^{r,\epsilon}$, and

$$\bar{\xi}_{r,0} \leq \kappa(\widehat{W}^r(0) \vee 1) \leq \kappa(\kappa_1 \vee 1) \leq \kappa(\kappa_1+\kappa_2+1).$$

Meanwhile, from (7.34), replacing $\epsilon$ by $\epsilon/\kappa(\kappa_1+\kappa_2+1)$, we have

$$|\widehat{Z}^r(t) - \zeta^\alpha \widehat{W}^r(t)| \leq \epsilon \qquad \text{for all } t \in [0, \bar{\xi}_{r,0}L/r].$$

Then from Proposition 9, with $r \geq 2\mathbf{J}/\epsilon_0$, we have

$$\int_0^{\bar{\xi}_{r,0}L/r} 1_{(\hat{W}^r(s)>1)}\, d\hat{Y}^r(s) = 0.$$

Now we assume that (7.39)–(7.42) hold up to $m$, and we shall show they also hold for $m+1$. First, (7.41) for $m$ directly implies (7.39) for $m+1$ because $\bar{\xi}_{r,m}L > 1$. Next, because (5.7) and (7.42) hold, and by Theorem 5.1 of Williams (1998b), we have

$$\mathrm{Osc}(\widehat{W}^r, [0,(m+\bar{\xi}_{r,m}L)/r]) \leq \mathrm{Osc}(\widehat{X}^r, [0,(m+\bar{\xi}_{r,m}L)/r]) + 1 \leq \kappa_2 + 1.$$

Note that, because $m \leq rT$ and (7.40), we have $(m+\bar{\xi}_{r,m}L)/r \leq T+L$ for $r$ large enough. It then follows that

$$\begin{aligned}
\bar{\xi}_{r,m+1} &\leq \kappa(\widehat{W}^r((m+1)/r) \vee 1) \\
&\leq \kappa(\widehat{W}^r(0) + \mathrm{Osc}(\widehat{W}^r, [0,(m+1)/r]) \vee 1) \\
&\leq \kappa(\widehat{W}^r(0) + \mathrm{Osc}(\widehat{W}^r, [0,(m+\bar{\xi}_{r,m}L)/r]) \vee 1) \\
&\leq \kappa(\kappa_1+\kappa_2+1).
\end{aligned}$$

This fact and (7.30) imply that

$$|\widehat{Z}^r(t) - \zeta^\alpha \widehat{W}^r(t)| \le \epsilon \qquad \text{for all } t \in [(m+1+\bar{\xi}_{r,m+1}\tau_0)/r, (m+1+\bar{\xi}_{r,m+1}L)/r].$$

We choose $L \ge 3N\tau_0 + 1$. Then Lemma 8 implies $\bar{\xi}_{r,m}L/r \ge (1+\bar{\xi}_{r,m+1}\tau_0)/r$. It follows that

$$|\widehat{Z}^r(t) - \zeta^\alpha \widehat{W}^r(t)| \le \epsilon \qquad \text{for all } t \in [0, (m+1+\bar{\xi}_{r,m+1}L)/r].$$

Then again Proposition 9 gives

$$\int_0^{(m+1+\bar{\xi}_{r,m+1}L)/r} 1_{(\widehat{W}^r(s)>1)}\, d\widehat{Y}^r(s) = 0.$$

Therefore, we can conclude that $\bar{\xi}_{r,m} \le \kappa(\kappa_1 + \kappa_2 + 1)$ for all $0 \le m \le rT$, which implies

$$\|\widehat{Z}^r(t) - \zeta^\alpha \widehat{W}^r(t)\|_T \le \epsilon \qquad \text{on } \mathcal{N}^{r,\epsilon}$$

for all large enough $r$. Theorem 5 follows because $\epsilon$ can be chosen arbitrarily small. □

**8. Heavy traffic limit theorem.** In this section we prove the heavy traffic limit theorem, Theorem 4. Our proof employs an invariance principle developed in Williams (1998b) for semimartingale reflecting Brownian motions (RBMs), specialized to the one-dimensional case.

Recall that the constant $\epsilon_0$ is defined in Proposition 9. For each $r$, define process $\delta^r = \{\delta^r(t) : t \ge 0\}$ as

$$\delta^r(t) = (|\widehat{Z}^r(t) - \zeta^\alpha \widehat{W}^r(t)| \vee (2\mathbf{J}/r))/\epsilon_0.$$

It follows from the state space collapse theorem, Theorem 5, that $\delta^r \to 0$ in probability as $r \to \infty$. Now, for each $r$, define processes $\gamma^r = \{\gamma^r(t) : t \ge 0\}$ and $\tilde{W}^r = \{\tilde{W}^r(t) : t \ge 0\}$ as

$$\gamma^r(t) = \widehat{W}^r(t) \wedge \delta^r(t) \quad \text{and} \quad \tilde{W}^r(t) = \widehat{W}^r(t) - \gamma^r(t).$$

Since $\delta^r \to 0$ in probability as $r \to \infty$, $\gamma^r \to 0$ in probability as $r \to \infty$. It is easy to see that $\tilde{W}^r(t) = (\widehat{W}^r(t) - \delta^r(t)) \vee 0 \ge 0$ for all $t \ge 0$. Now we claim that

$$\int_0^\infty \tilde{W}^r(s)\, d\widehat{Y}^r(s) = 0. \tag{8.1}$$

To see this, it is enough to show that for any $0 \le t_1 < t_2$, $\widehat{Y}^r(t_2) = \hat{Y}^r(t_1)$ whenever $\tilde{W}^r(t) > 0$ for all $t \in [t_1, t_2]$; see, for example, Lemma 3.1.2 of Dai and Williams (2003). Suppose that $\tilde{W}^r(t) > 0$ for all $t \in [t_1, t_2]$. Then $\widehat{W}^r(t) > \delta^r(t)$ for all $t \in [t_1, t_2]$. The latter condition implies that $\widehat{W}^r(t) \ge$

$2\mathbf{J}/(r\epsilon_0)$ and $|\widehat{Z}^r(t)/\widehat{W}^r(t) - \zeta^\alpha| \le \epsilon_0$ for all $t \in [t_1, t_2]$. It follows from Proposition 9 that $\widehat{Y}^r(t_2) = \widehat{Y}^r(t_1)$, thus proving (8.1).

By Theorem 1, the maximum pressure policy is asymptotically efficient, and thus, by Lemma 3, $\widehat{X}^r$ converges in distribution to the Brownian motion $X^*$ as $r \to \infty$.

Therefore, we conclude that the processes $(\widehat{W}^r, \widehat{X}^r, \widehat{Y}^r)$ satisfy the following conditions:

(i) $\widehat{W}^r = \widehat{X}^r + \widehat{Y}^r$,

(ii) $\widehat{W}^r = \tilde{W}^r + \gamma^r$, where $\tilde{W}^r(t) \ge 0$ for all $t \ge 0$, and $\gamma^r \to 0$ in probability as $r \to \infty$,

(iii) $\widehat{X}^r$ converges to the Brownian motion $X^*$ in distribution as $r \to \infty$.

(iv) With probability 1,

- (a) $\widehat{Y}^r(0) = 0$,
- (b) $\widehat{Y}^r$ is nondecreasing,
- (c) $\int_0^\infty \tilde{W}^r(s)\, d\widehat{Y}^r(s) = 0$.

Since we have a one-dimensional case, condition (II) of Proposition 4.2 of Williams (1998b) is satisfied, and consequently, condition (vi) of Theorem 4.1 of Williams (1998b) holds. It then follows from Theorem 4.1 of Williams (1998b) that $\widehat{W}^r$ converges in distribution to the reflecting Brownian motion $W^* = \psi(X^*)$. The convergence of $\widehat{Z}^r$ to $\zeta^\alpha W^*$ follows from the state space collapse result (5.16). Thus, we have completed the proof of Theorem 4. □

**9. Linear holding cost.** In this section we discuss the asymptotic optimality in terms of a linear holding cost structure for the sequence of stochastic processing networks. For each network in the sequence, we assume that there is a linear holding cost incurred at rate $c_i > 0$ for each job in buffer $i$. Let $c$ be the corresponding vector. For the $r$th network, the total holding cost rate at time $t$ is

$$C^r(t) = c \cdot Z^r(t). \tag{9.1}$$

Define the diffusion-scaled linear holding cost rate processes of the $r$th network $\widehat{C}^r = \{\widehat{C}^r(t), t \ge 0\}$ via

$$\widehat{C}^r(t) = C^r(r^2 t)/r.$$

Clearly, $\widehat{C}^r(t) = c \cdot \widehat{Z}^r(t)$ for $t \ge 0$.

DEFINITION 7. Consider a sequence of stochastic processing networks indexed by $r$. For a given $\varepsilon > 0$, an asymptotically efficient policy $\pi^*$ is said to be *asymptotically $\varepsilon$-optimal for the linear holding cost* if for any $t>0, \eta>0$, and any asymptotically efficient policy $\pi$,

$$\limsup_{r\to\infty} \mathbb{P}(\widehat{C}^r_{\pi^*}(t)>\eta) \le \liminf_{r\to\infty} \mathbb{P}((1+\varepsilon)\widehat{C}^r_{\pi}(t)>\eta), \tag{9.2}$$

where $\widehat{C}^r_{\pi^*}(t)$ and $\widehat{C}^r_{\pi}(t)$ are the diffusion-scaled total holding cost rates at time $t$ under policies $\pi^*$ and $\pi$, respectively.

THEOREM 7. *Consider a sequence of stochastic processing networks where Assumptions 1–4 hold and the limit network satisfies the EAA condition. For any given $\varepsilon>0$, there exists a maximum pressure policy $\pi^*$ that is asymptotically $\varepsilon$-optimal for the linear holding cost.*

PROOF. Firs, we note that under any policy, for all $r$,

$$\widehat{C}^r(t) = c\cdot \widehat{Z}^r(t) \ge \left(\min_{i\in\mathcal{I}} c_i/y^r_i\right)\widehat{W}^r(t).$$

Let $\iota \in \arg\min_{i\in\mathcal{I}} c_i/y^*_i$. Then for any asymptotically efficient policy $\pi$,

$$\liminf_{r\to\infty}\mathbb{P}(\widehat{C}^r_{\pi}(t)>\eta) \ge \liminf_{r\to\infty}\mathbb{P}((c_\iota/y^r_\iota)\widehat{W}^r_\pi(t)>\eta) \ge \mathbb{P}((c_\iota/y^*_\iota)W^*(t)>\eta).$$

The second inequality follows from (5.10) and the fact that $y^r\to y^*$ as $r\to\infty$.

Now consider a maximum pressure policy $\pi^*$ with parameter $\alpha$ given by

$$\alpha_i = \begin{cases} c_i y^*_i\varepsilon/(\mathbf{I}|c\times y^*|), & i=\iota,\\ 1, & \text{otherwise.}\end{cases}$$

Because $\sum_{i'}(y^*_{i'})^2/\alpha_{i'} \ge (y^*_\iota)^2/\alpha_\iota = y^*_\iota \mathbf{I}|c\times y^*|/c_\iota\varepsilon$, we have $c_i\zeta^\alpha_i \le c_\iota\varepsilon/(y^*_\iota\mathbf{I})$ for all $i\ne\iota$ and $c_\iota\zeta^\alpha_\iota \le c_\iota/y^*_\iota$.

From Theorem 4, we have

$$\widehat{C}^r_{\pi^*} \quad\Rightarrow\quad c\cdot Z^{*,\alpha} \qquad \text{as } r\to\infty.$$

Because

$$\sum_i c_i Z^{*,\alpha}_i(t) = \sum_i c_i\zeta^\alpha_i W^*(t) \le (c_\iota/y^*_\iota)W^*(t)(1+\varepsilon),$$

we have

$$\lim_{r\to\infty}\mathbb{P}(\widehat{C}^r_{\pi^*}(t)>\eta) \le \mathbb{P}((c_\iota/y_\iota)W^*(t)(1+\varepsilon)>\eta).$$

Then (9.2) follows. □

Theorem 7 says that, for any $\varepsilon>0$, one can always find a maximum pressure policy such that, for any time $t$, the total holding cost rate under the maximum pressure policy is asymptotically dominated by $(1+\varepsilon)$ times the total holding cost rate under any other efficient policy in the sense of stochastic ordering. For a given $\epsilon$, the parameter $\alpha$ that is used to define the asymptotically $\varepsilon$-optimal maximum pressure policy in Theorem 7 depends on the first-order network data $R$ and $A$.

## APPENDIX A: AN EQUIVALENT DUAL FORMULATION FOR THE STATIC PLANNING PROBLEM

In this section we describe an equivalent dual formulation for the static planning problem (3.1)–(3.5). The equivalent dual formulation will be useful for our proofs in Appendix B. Throughout this section, we will consider an arbitrary stochastic processing network, so the results developed here can be applied to each of the $r$th network in the network sequence that we discussed in Section 4 and they can also be applied to the limit network. The main result of this section is the following proposition.

PROPOSITION A.1. *Suppose $\rho$ is the optimal objective value to the static planning problem (3.1)–(3.5). Then $(y^*, z^*)$ is optimal to the dual problem (3.6)–(3.10) if and only if $y^*$ satisfies*

$$\max_{a\in\mathcal{A}} \sum_{i\in\mathcal{I},j\in\mathcal{J}_I} y_i^* R_{ij} a_j = -\rho \tag{A.1}$$

*and*

$$\max_{a\in\mathcal{A}} \sum_{i\in\mathcal{I},j\in\mathcal{J}_S} y_i^* R_{ij} a_j = 1, \tag{A.2}$$

*and $\{z_k^*, k\in\mathcal{K}_I\}$ and $\{z_k^*, k\in\mathcal{K}_S\}$ are, respectively, optimal solutions to*

$$\text{(A.3)}\quad \min \quad -\sum_{k\in\mathcal{K}_I} z_k$$

$$\text{(A.4)}\quad \text{s.t.} \quad \sum_{i\in\mathcal{I}} y_i^* R_{ij} \le -\sum_{k\in\mathcal{K}_I} z_k A_{kj} \qquad \textit{for each input activity } j;$$

*and*

$$\text{(A.5)}\quad \min \quad \sum_{k\in\mathcal{K}_S} z_k$$

$$\text{(A.6)}\quad \text{s.t.} \quad \sum_{i\in\mathcal{I}} y_i^* R_{ij} \le \sum_{k\in\mathcal{K}_S} z_k A_{kj} \qquad \textit{for each service activity } j,$$

$$\text{(A.7)}\quad z_k \ge 0 \qquad \textit{for each service processor.}$$

PROOF. For the "only if" part, we assume $(x^*, \rho)$ and $(y^*, z^*)$ are an optimal dual pair for the static planning problem and its dual problem. We will show that $y^*$ satisfies (A.1)–(A.2) and $z^*$ is optimal to (A.3)–(A.4) and (A.5)–(A.7). We first show that

$$\max_{a\in\mathcal{A}} \sum_{i\in\mathcal{I},j\in\mathcal{J}_I} y_i^* R_{ij} a_j \ge -\rho, \tag{A.8}$$

and

$$\max_{a \in \mathcal{A}} \sum_{i \in \mathcal{I}, j \in \mathcal{J}_S} y_i^* R_{ij} a_j \geq 1. \tag{A.9}$$

For this, we construct a feasible allocation $a^*$ with

$$a_j^* = \begin{cases} x_j^*, & j \in \mathcal{J}_I, \\ x_j^*/\rho, & j \in \mathcal{J}_S. \end{cases}$$

By the complementary slackness on the constraints (3.7) and (3.8), we have

$$\sum_{i \in \mathcal{I}, j \in \mathcal{J}_I} y_i^* R_{ij} x_j^* = - \sum_{k \in \mathcal{K}_I, j \in \mathcal{J}_I} z_k^* A_{kj} x_j^* \tag{A.10}$$

and

$$\sum_{i \in \mathcal{I}, j \in \mathcal{J}_S} y_i^* R_{ij} x_j^* = \sum_{k \in \mathcal{K}_S, j \in \mathcal{J}_S} z_k^* A_{kj} x_j^*. \tag{A.11}$$

By the complementary slackness on the constraints (3.3) and (3.4), we have

$$\sum_{j \in \mathcal{J}_I, k \in \mathcal{K}_I} z_k^* A_{kj} x_j^* = \sum_{k \in \mathcal{K}_I} z_k^* = \rho \tag{A.12}$$

and

$$\sum_{j \in \mathcal{J}_S, k \in \mathcal{K}_S} z_k^* A_{kj} x_j^* = \rho \sum_{k \in \mathcal{K}_S} z_k^* = \rho. \tag{A.13}$$

The last equality in (A.12) is from the strong duality theorem; the optimal objective value of the dual problem equals the optimal objective value of the primal problem. Readers are referred to Section 4.2 of Luenberger (1984) for a formal description of the strong duality theorem. The last equality in (A.13) follows from (3.9). Then from the definition of $a^*$ and (A.10)–(A.13), it immediately follows that

$$\sum_{i \in \mathcal{I}, j \in \mathcal{J}_I} y_i^* R_{ij} a_j^* = \sum_{i \in \mathcal{I}, j \in \mathcal{J}_I} y_i^* R_{ij} x_j^* = -\rho$$

and

$$\sum_{i \in \mathcal{I}, j \in \mathcal{J}_S} y_i^* R_{ij} a_j^* = \sum_{i \in \mathcal{I}, j \in \mathcal{J}_S} y_i^* R_{ij} x_j^* / \rho = 1,$$

which imply (A.8) and (A.9) because $a^* \in \mathcal{A}$.

Next we shall show that

$$\max_{a \in \mathcal{A}} \sum_{i \in \mathcal{I}, j \in \mathcal{J}_I} y_i^* R_{ij} a_j \leq -\rho$$

and

$$\max_{a\in\mathcal{A}} \sum_{i\in\mathcal{I},j\in\mathcal{J}_S} y_i^* R_{ij} a_j \le 1.$$

For any $a \in \mathcal{A}$, we have

$$\sum_{i\in\mathcal{I},j\in\mathcal{J}_I} y_i^* R_{ij} a_j \le -\sum_{k\in\mathcal{K}_I,j\in\mathcal{J}_I} z_k^* A_{kj} a_j = -\sum_{k\in\mathcal{K}_I} z_k^* = -\rho.$$

The first inequality above follows from (3.7), and the nonnegativity of $a$; the second inequality holds since $a \in \mathcal{A}$ and, therefore, $\sum_{j\in\mathcal{J}_I} A_{kj} a_j = 1$ for each $k \in \mathcal{K}_I$; the third is due to the strong duality theorem. Similarly, we have

$$\sum_{i\in\mathcal{I},j\in\mathcal{J}_S} y_i^* R_{ij} a_j \le \sum_{k\in\mathcal{K}_S,j\in\mathcal{J}_S} z_k^* A_{kj} a_j \le \sum_{k\in\mathcal{K}_S} z_k^* = 1.$$

The first inequality above follows from (3.8); the second follows from $\sum_{j\in\mathcal{J}_S} A_{kj} a_j \le 1$ and $z_k \ge 0$ for each $k \in \mathcal{K}_S$; and the third is due to (3.9).

Hence, $y^*$ satisfies (A.1) and (A.2). To see the $z^*$ is an optimal solution to (A.3)–(A.4) and (A.5)–(A.7), we consider the following problems:

$$\text{(A.14)}\qquad \begin{aligned} \max \quad & \sum_{i\in\mathcal{I},j\in\mathcal{J}_I} y_i^* R_{ij} a_j \\ \text{s.t.} \quad & \sum_{j\in\mathcal{J}_I} A_{kj} a_j = 1 \qquad \text{for each input processor } k, \\ & a_j \ge 0, \qquad j \in \mathcal{J}_I \end{aligned}$$

and

$$\text{(A.15)}\qquad \begin{aligned} \max \quad & \sum_{i\in\mathcal{I},j\in\mathcal{J}_S} y_i^* R_{ij} a_j \\ \text{s.t.} \quad & \sum_{j\in\mathcal{J}_S} A_{kj} a_j \le 1 \qquad \text{for each service processor } k, \\ & a_j \ge 0, \qquad j \in \mathcal{J}_S. \end{aligned}$$

It is easy to see that the above two problems are equivalent to the left-hand side of (A.1) and (A.2). Furthermore, they are the dual problems of (A.3)–(A.4) and (A.5)–(A.7). This implies that the optimal objective values to (A.3)–(A.4) and (A.5)–(A.7) are $-\rho$ and 1 respectively. Because $(y^*, z^*)$ is an optimal solution to (3.6)–(3.10), (A.4) and (A.6) are satisfied by $z^*$, $\sum_{k\in\mathcal{K}_I} z_k^* = \rho$, and $\sum_{k\in\mathcal{K}_S} z_k^* = 1$. This implies that $z^*$ is feasible to (A.3)–(A.4) and (A.5)–(A.7) with respective objective values $-\rho$ and 1. Therefore, $z^*$ is optimal to (A.3)–(A.4) and (A.5)–(A.7).

For the "if" part, let $(y^*, z^*)$ be such that $y^*$ satisfies (A.1) and (A.2) and $z^*$ is optimal to (A.3)–(A.4) and (A.5)–(A.7). Because (A.3)–(A.4) and

(A.5)–(A.7) are dual problems of the equivalent formulation of the left-hand side of (A.1) and (A.2),

$$\sum_{k\in\mathcal{K}_I} z_k^* = \rho \tag{A.16}$$

and

$$\sum_{k\in\mathcal{K}_S} z_k^* = 1. \tag{A.17}$$

The fact that $z^*$ is feasible to (A.3)–(A.4) and (A.5)–(A.7), together with (A.17), implies that $(y^*, z^*)$ is feasible to the dual problem (3.6)–(3.10). Furthermore, the corresponding objective value is $\rho$ because of (A.16). This implies that $(y^*, z^*)$ is optimal to (3.6)–(3.10). □

Proposition A.1 immediately leads to the following corollary because the problem $\max_{a\in\mathcal{A}} y^*\cdot Ra$ can be decomposed into problems (A.14) and (A.15).

COROLLARY A.1. *Suppose $(y^*, z^*)$ is the unique optimal solution to the dual problem (3.6)–(3.10) with objective value $\rho$. Then $y^*$ is the unique* **I**-*dimensional vector that satisfies*

$$\max_{a\in\mathcal{A}} \sum_{i\in\mathcal{I}, j\in\mathcal{J}_S} y_i^* R_{ij} a_j = 1 \tag{A.18}$$

*and*

$$\max_{a\in\mathcal{A}} y^*\cdot Ra = 1-\rho. \tag{A.19}$$

## APPENDIX B: PROOFS OF LEMMAS AND THEOREM 1

In this section we provide the proofs for Lemmas 1, 2, 5, 6 and 10, and Theorem 1.

PROOF OF LEMMA 2. From Corollary A.1 and Assumption 2, we have

$$\max_{a\in\mathcal{A}} y^r\cdot Ra = 1-\rho^r.$$

On the other hand, $T^r(t)/t \in \mathcal{A}$. Hence, $(1-\rho^r)t \ge y^r\cdot RT^r(t)$, which implies $Y^r$ is nonnegative. Similarly, for any $0 \le t_1 < t_2$, $(T^r(t_2)-T^r(t_1))/(t_2-t_1) \in \mathcal{A}$, and we have

$$y^r\cdot R(T^r(t_2)-T^r(t_1))/(t_2-t_1) \le (1-\rho^r).$$

It follows that $Y^r(t_2)-Y^r(t_1) \ge 0$. □

PROOF OF LEMMA 5. First, if the statement is not true, then we can find a $v_0 \in V^o$ such that $\kappa v_0 \notin V$ for all $0 < \kappa \leq 1$ because of the convexity of $V$. Denote $V_0 = \{\kappa v_0, 0 < \kappa < 1\}$. Because any $v$ in $V_0$ is not in $V$, $V_0 \cap V = \varnothing$. It is easy to see that $V_0$ is relatively open and convex. Therefore, there exists a hyperplane separating $V$ and $V_0$ [cf. Rudin (1991), Theorem 3.4]. In other words, there exists a vector $\hat{y}$ and a constant $b$ such that $\hat{y} \cdot v \leq b$ for all $v \in V$ and $\hat{y} \cdot v > b$ for all $v \in V_0$. We notice that $b$ must be zero. To see this, first we have $b \geq 0$ because the origin is in $V$. Moreover, for any $\epsilon > 0$, we can choose $\kappa$ arbitrarily small such that $\kappa \hat{y} \cdot v_0 < \epsilon$. Because $\kappa v_0 \in V_0$, we have $b < \kappa \hat{y} \cdot v_0 < \epsilon$. This implies $b = 0$, therefore, the origin is in the separating hyperplane and $\max_{a \in \mathcal{A}} \hat{y} \cdot Ra = 0$. Obviously, $\hat{y} \neq y^*$ because $y^* \cdot v = 0 < \hat{y} \cdot v$ for $v \in V_0$. Since $\max_{a \in \mathcal{A}} \sum_{i \in \mathcal{I}, j \in \mathcal{J}_S} \hat{y}_i R_{ij} a_j > 0$ for any vector $\hat{y}$, we consider two cases:

Case 1: $\max_{a \in \mathcal{A}} \sum_{i \in \mathcal{I}, j \in \mathcal{J}_S} \hat{y}_i R_{ij} a_j > 0$. For this case, without loss of generality, we select $\hat{y}$ such that $\max_{a \in \mathcal{A}} \sum_{i \in \mathcal{I}, j \in \mathcal{J}_S} \hat{y}_i R_{ij} a_j = 1$. Then $\hat{y}$ satisfies both (A.18) and (A.19) with $\rho = 1$. On the other hand, from Corollary A.1, $y^*$ is the unique vector that satisfies both (A.18) and (A.19). This is a contradiction.

Case 2: $\max_{a \in \mathcal{A}} \sum_{i \in \mathcal{I}, j \in \mathcal{J}_S} \hat{y}_i R_{ij} a_j = 0$. For this case, we will show that $y^* + \hat{y}$ satisfies both (A.18) and (A.19) with $\rho = 1$, thus yielding a contradiction. For this, first observe that $x^* \in \arg\max_{a \in \mathcal{A}} y^* \cdot Ra$. This implies that

$$x^* \in \arg\max_{a \in \mathcal{A}} \sum_{i \in \mathcal{I}, j \in \mathcal{J}_S} y_i^* R_{ij} a_j.$$

To see this, suppose there is an allocation $\tilde{a} \in \mathcal{A}$ such that

$$\sum_{i \in \mathcal{I}, j \in \mathcal{J}_S} y_i^* R_{ij} \tilde{a}_j > \sum_{i \in \mathcal{I}, j \in \mathcal{J}_S} y_i^* R_{ij} x_j^*,$$

then we can define another allocation $\hat{a}$ by

$$\hat{a}_j = \begin{cases} x_j^*, & j \in \mathcal{J}_I, \\ \tilde{a}_j, & j \in \mathcal{J}_S, \end{cases}$$

so that $y^* \cdot R\hat{a} > y^* \cdot Rx^*$. This is certainly not true.

We also observe that $x^* \in \arg\max_{a \in \mathcal{A}} \hat{y} \cdot Ra$ because $\max_{a \in \mathcal{A}} \hat{y} \cdot Ra = 0$ and $x^*$ satisfies (3.2). This again implies that

$$x^* \in \arg\max_{a \in \mathcal{A}} \sum_{i \in \mathcal{I}, j \in \mathcal{J}_S} \hat{y}_i R_{ij} a_j.$$

Thus,

$$\max_{a \in \mathcal{A}} (y^* + \hat{y}) \cdot Ra = (y^* + \hat{y}) \cdot Rx^* = 0$$

and

$$\max_{a\in\mathcal{A}} \sum_{i\in\mathcal{I},j\in\mathcal{J}_S} (y_i^* + \hat{y}_i)R_{ij}a_j = \sum_{i\in\mathcal{I},j\in\mathcal{J}_S} (y_i^* + \hat{y}_i)R_{ij}x_j^* = 1. \qquad \square$$

PROOF OF LEMMA 6. It is natural to work in a more general setting. Consider an i.i.d. sequence of nonnegative random variables $\{v_\ell, \ell = 1, 2, \ldots\}$ with mean $1/\mu_v$. Assume $v_\ell$ have finite $2 + \epsilon_v$ moments for some $\epsilon_v > 0$. That is, there exists some $\hat{\sigma} < \infty$ such that $\mathbb{E}(v_\ell^{2+\epsilon_v}) = \hat{\sigma}$. Let $V(n) = \sum_{\ell=1}^n v_\ell, n \in \mathbb{Z}^+$. Define the renewal process associated with $V(n)$ as $G(t) = \max\{n : V(n) \le t\}$. Let $v^{r,T,\max} = \max\{v_\ell : 1 \le \ell \le G(r^2T) + 1\}$.

It immediately follows from Lemma 3.3 of Iglehart and Whitt (1970) that

$$v^{r,T,\max}/r \to 0 \qquad \text{with probability 1.} \tag{B.1}$$

We now show that for any fixed $\epsilon > 0$ and large enough $n$,

$$\mathbb{P}(\|V(\ell) - \ell/\mu_v\|_n \ge \epsilon n) \le \epsilon/n. \tag{B.2}$$

Because $v_\ell$ have finite $2 + \epsilon_v$ moments, one gets

$$\mathbb{E}(|V(\ell) - \ell/\mu_v|^{2+\epsilon_v}) \le \kappa_v \ell^{1+\epsilon_v/2} \qquad \text{for all } \ell \le n,$$

where $\kappa_v$ is some constant that depends just on $\hat{\sigma}$ and $\epsilon_v$ [cf. Ata and Kumar (2005), Lemma 8]. Then from Chebyshev's inequality, we have, for each $\ell \le n$,

$$\mathbb{P}(|V(\ell) - \ell/\mu_v| \ge \epsilon n) \le \kappa_v \ell^{1+\epsilon_v/2}/(\epsilon n)^{2+\epsilon_v} \le \kappa_v/(\epsilon^{2+\epsilon_v} n^{1+\epsilon_v/2}).$$

Choosing $n$ large enough,

$$\mathbb{P}(|V(\ell) - \ell/\mu_v| \ge \epsilon n) \le \epsilon/n.$$

Let

$$\tau = \min\{\ell : |V(\ell) - \ell/\mu_v| \ge n\epsilon\}.$$

Then, restarting the process at time $\tau$,

$$\begin{aligned}&\mathbb{P}(|V(n) - n/\mu_v| \le \epsilon n/2 \mid \tau \le n)\\ &\qquad \le \mathbb{P}(|V(n) - V(\tau) - (n - \tau)/\mu_v| \ge \epsilon n/2 \mid \tau \le n) \le \epsilon/2n.\end{aligned}$$

On the other hand,

$$\mathbb{P}(|V(n) - n/\mu_v| \le \epsilon n/2) \ge 1 - \epsilon/2n.$$

It then follows that $\mathbb{P}(\tau \le n) \le \epsilon/(2n - \epsilon) \le \epsilon/n$. This implies the result.

Finally, the inversion of (B.2) gives

$$\mathbb{P}(\|G(s) - \mu_v s\|_t \ge \epsilon t) \le \epsilon/t \qquad \text{for large enough } t. \tag{B.3}$$

Readers are referred to Proposition 4.3 of Bramson (1998) for such an inversion.

Applying (B.1) and (B.3) to the utilized service times $u_j(\ell)$, one gets (7.4) and (7.6) in Lemma 6. Using (B.2) for each component of the routing vector $\phi_i^j(\ell)$ yields (7.5). □

Proof of Lemma 1. We first show that $\{(y^r, z^r)\}$ is bounded. Since $z^r \geq 0$, $\sum_{k\in\mathcal{K}_S} z_k^r = 1$, and $\sum_{k\in\mathcal{K}_I} z_k^r = \rho^r$ for all $r$, we have $|z^r| \leq 1$. To show $\{y^r\}$ is bounded, we consider the following primal-dual pair:

$$\text{minimize} \quad \rho \tag{B.4}$$

$$\text{subject to} \quad Rx \geq be, \tag{B.5}$$

$$\sum_{j\in\mathcal{J}} A_{kj} x_j = 1 \quad \text{for each input processor } k, \tag{B.6}$$

$$\sum_{j\in\mathcal{J}} A_{kj} x_j \leq \rho \quad \text{for each service processor } k, \tag{B.7}$$

$$x \geq 0, \tag{B.8}$$

and

$$\text{maximize} \quad \sum_{k\in\mathcal{K}_I} z_k + b \sum_{i\in\mathcal{I}} y_i, \tag{B.9}$$

$$\text{subject to} \quad \sum_{i\in\mathcal{I}} y_i R_{ij} \leq - \sum_{k\in\mathcal{K}_I} z_k A_{kj} \quad \text{for each input activity } j, \tag{B.10}$$

$$\sum_{i\in\mathcal{I}} y_i R_{ij} \leq \sum_{k\in\mathcal{K}_S} z_k A_{kj} \quad \text{for each service activity } j, \tag{B.11}$$

$$\sum_{k\in\mathcal{K}_S} z_k = 1, \tag{B.12}$$

$$y \geq 0; \quad \text{and} \quad z_k \geq 0 \quad \text{for each service processor } k. \tag{B.13}$$

The dual LP (B.9)–(B.13) is obtained by perturbing the objective function coefficients of the dual static planning problem (3.6)–(3.10). Because the dual static planning problem (3.6)–(3.10) has a unique optimal solution, for sufficiently small $b > 0$, the optimal solution of the dual LP (B.9)–(B.13) equals $(y^*, z^*)$ [cf. Mangasarian (1979), Theorem 1]. Therefore, the primal problem (B.4)–(B.8) has an optimal solution $(\hat{\rho}, \hat{x})$. Now choose $r$ large enough such that $|R^r\hat{x} - R\hat{x}| < be/2$. Then $R^r\hat{x} \geq be/2$. Consider the problem (B.4)–(B.8) with $b$ replaced by $b/2$ and $R$ replaced by $R^r$. Because $(\hat{\rho}, \hat{x})$ is a feasible solution, the optimal objective value $\tilde{\rho}^r \leq \hat{\rho}$. The corresponding dual problem of this new LP is the dual problem (B.9)–(B.13) with $b$ in the objective function coefficients replaced by $b/2$ and $R$ replaced by $R^r$. The optimal objective value of this new dual LP equals $\tilde{\rho}^r \le \hat{\rho}$. The new dual LP has the exact same constraints as the dual static planning problem (3.6)–(3.10) for the $r$th network. Thus, $(y^r, z^r)$ is a feasible solution to the new dual LP. It then follows that

$$\sum_{k\in\mathcal{K}_I} z_k^r + b/2 \sum_{i\in\mathcal{I}} y_i^r \le \tilde{\rho}^r \le \hat{\rho}.$$

This implies that $\sum_i y_i^r \le 2\hat{\rho}/b$ for large enough $r$, so $\{y^r\}$ is bounded.

Then we only need to show that every convergent subsequence of $\{(y^r, z^r)\}$ converges to $(y^*, z^*)$. Let $(\hat{y}, \hat{z})$ be a limit point of any subsequence $\{(y^{r_n}, z^{r_n})\}$. We will verify that $(\hat{y}, \hat{z})$ is an optimal solution to the dual static planning problem (3.6)–(3.10) of the limiting network. First, we show that they are feasible. Since $\{(y^{r_n}, z^{r_n})\} \to (\hat{y}, \hat{z})$ as $n \to \infty$, for any $\epsilon > 0$, for large enough $n$, $|\hat{y} - y^{r_n}| < \epsilon$, $|\hat{z} - z^{r_n}| < \epsilon$ and $|R - R^{r_n}| < \epsilon$. For each input activity $j \in \mathcal{J}_I$,

$$\begin{aligned}\sum_{i\in\mathcal{I}} \hat{y}_i R_{ij} &\le \sum_{i\in\mathcal{I}} y_i^{r_n} R_{ij}^{r_n} + \mathbf{I}\epsilon\Big(|\hat{y}| + \sup_r |R^r|\Big)\\ &\le -\sum_{k\in\mathcal{K}_I} A_{kj} z_k^{r_n} + \mathbf{I}\epsilon\Big(|\hat{y}| + \sup_r |R^r|\Big)\\ &\le -\sum_{k\in\mathcal{K}_I} A_{kj} \hat{z}_k + \mathbf{I}\epsilon\Big(|\hat{y}| + \sup_r |R^r|\Big) + \mathbf{K}\epsilon.\end{aligned}$$

Since $\epsilon$ can be arbitrarily small, we have

$$\sum_{i\in\mathcal{I}} \hat{y}_i R_{ij} \le -\sum_{k\in\mathcal{K}_I} A_{kj} \hat{z}_k \qquad \text{for each input activity } j \in \mathcal{J}_I.$$

Similarly, one can verify that

$$\sum_{i\in\mathcal{I}} \hat{y}_i R_{ij} \le \sum_{k\in\mathcal{K}_S} A_{kj} \hat{z}_k \qquad \text{for each service activity } j \in \mathcal{J}_S$$

and

$$\sum_{k\in\mathcal{K}_S} \hat{z}_k = 1.$$

Furthermore, because $(y^{r_n}, z^{r_n})$ are optimal solutions, $\sum_{k\in\mathcal{K}_I} z_k^{r_n} = \rho^{r_n}$. Again, we can show that

$$\sum_{k\in\mathcal{K}_I} \hat{z}_k = 1.$$

Therefore, $(\hat{y}, \hat{z})$ is an optimal solution to the dual problem (3.6)–(3.10) of the limit network. Then by the uniqueness of the optimal solution, we conclude that $(\hat{y},\hat{z})=(y^*,z^*)$. Since the subsequence is arbitrary, we have $(y^r,z^r)\to(y^*,z^*)$ as $r\to\infty$. □

PROOF OF LEMMA 10. Similar to the proof of Lemma 1 above, we can prove that $x^r\to x^*$ as $r\to\infty$. From the strict complementary theorem [cf. Wright (1997)], every pair of primal and dual LPs has a strict complementary optimal solution if they have optimal solutions. Hence, we have the following relations: for the limit network,

$$\sum_{j\in\mathcal{J}_S} A_{kj}x_j^*=1 \text{ is equivalent to } z_k^*>0 \qquad \text{for all } k\in\mathcal{K}_S; \tag{B.14}$$

$$\sum_{i\in\mathcal{I}} y_i^*R_{ij}=\sum_{k\in\mathcal{K}_S} A_{kj}z_k^* \text{ is equivalent to } x_j^*>0 \qquad \text{for all } j\in\mathcal{J}_S; \tag{B.15}$$

$$\sum_{i\in\mathcal{I}} y_i^*R_{ij}=-\sum_{k\in\mathcal{K}_I} A_{kj}z_k^* \text{ is equivalent to } x_j^*>0 \qquad \text{for all } j\in\mathcal{J}_I; \tag{B.16}$$

and for each $r$,

$$z_k^r>0 \quad \text{implies} \quad \sum_{j\in\mathcal{J}_S} A_{kj}x_j^r=\rho^r \qquad \text{for all } k\in\mathcal{K}_S; \tag{B.17}$$

$$x_j^r>0 \quad \text{implies} \quad \sum_{i\in\mathcal{I}} y_i^rR_{ij}^r=\sum_{k\in\mathcal{K}_S} A_{kj}z_k^r \qquad \text{for all } j\in\mathcal{J}_S; \tag{B.18}$$

$$x_j^r>0 \quad \text{implies} \quad \sum_{i\in\mathcal{I}} y_i^rR_{ij}^r=-\sum_{k\in\mathcal{K}_I} A_{kj}z_k^r \qquad \text{for all } j\in\mathcal{J}_I. \tag{B.19}$$

Since $x^r\to x^*$ as $r\to\infty$, we have for large enough $r$,

$$x_j^*>0 \text{ implies } x_j^r>0 \qquad \text{for all } j\in\mathcal{J}. \tag{B.20}$$

This, together with (B.15) and (B.18), implies that, for large enough $r$ and each service activity $j\in\mathcal{J}_S$,

$$\sum_{i\in\mathcal{I}} y_i^*R_{ij}=\sum_{k\in\mathcal{K}_S} A_{kj}z_k^* \text{ implies } \sum_{i\in\mathcal{I}} y_i^rR_{ij}^r=\sum_{k\in\mathcal{K}_S} A_{kj}z_k^r. \tag{B.21}$$

Similarly, it follows from (B.16), (B.19) and (B.20) that, for large enough $r$ and each input activity $j\in\mathcal{J}_I$,

$$\sum_{i\in\mathcal{I}} y_i^*R_{ij}=-\sum_{k\in\mathcal{K}_I} A_{kj}z_k^* \text{ implies } \sum_{i\in\mathcal{I}} y_i^rR_{ij}^r=-\sum_{k\in\mathcal{K}_I} A_{kj}z_k^r. \tag{B.22}$$

Suppose $z_k^*=0$ for some $k\in\mathcal{K}_S$, then $\sum_{j\in\mathcal{J}_S} A_{kj}x_j^*<1$. There exists an $\epsilon>0$ such that $\sum_{j\in\mathcal{J}_S} A_{kj}x_j^*=1-\epsilon$. For large enough $r$, we have $\sum_{j\in\mathcal{J}_S} A_{kj}x_j^r\le\sum_{j\in\mathcal{J}_S} A_{kj}x_j^*+\epsilon/2\le 1-\epsilon/2$ because $x^r\to x^*$ as $r\to\infty$. This implies $z_k^r=0$ for large enough $r$. Therefore, for large enough $r$,

$$z_k^*=0 \text{ implies } z_k^r=0 \qquad \text{for all } k\in\mathcal{K}_S. \tag{B.23}$$

Because $(y^*, z^*)$ is the optimal solution to the dual static planning problem (3.6)–(3.10), we have for each $a \in \mathcal{E}$,

$$\sum_{j\in\mathcal{J}}\left(\sum_{i\in\mathcal{I}} y_i^* R_{ij}\right) a_j \le \sum_{j\in\mathcal{J}_S}\left(\sum_{k\in\mathcal{K}_S} A_{kj} z_k^*\right) a_j - \sum_{j\in\mathcal{J}_I}\left(\sum_{k\in\mathcal{K}_I} A_{kj} z_k^*\right) a_j \le 0. \tag{B.24}$$

For the second inequality, we use the fact that $\sum_{j\in\mathcal{J}_S} A_{kj} a_j \le 1$ for all $k \in \mathcal{K}_S$, $\sum_{j\in\mathcal{J}_I} A_{kj} a_j = 1$ for all $k \in \mathcal{K}_I$, and $\sum_{k\in\mathcal{K}_S} z_k^* = \sum_{k\in\mathcal{K}_I} z_k^* = 1$. Since $a^* \in \mathcal{E}^*$, $y^* \cdot Ra^* = \max_{a\in\mathcal{E}} y^* \cdot Ra = 0$. It follows that both inequalities in (B.24) are equalities for $a^*$. Therefore,

$$\sum_{i\in\mathcal{I}} y_i^* R_{ij} = \sum_{k\in\mathcal{K}_S} A_{kj} z_k^* \qquad \text{for all } j \in \mathcal{J}_S \text{ with } a_j^* > 0, \tag{B.25}$$

$$\sum_{i\in\mathcal{I}} y_i^* R_{ij} = -\sum_{k\in\mathcal{K}_I} A_{kj} z_k^* \qquad \text{for all } j \in \mathcal{J}_I \text{ with } a_j^* > 0, \tag{B.26}$$

$$\sum_{j\in\mathcal{J}_S} A_{kj} a_j^* = 1 \qquad \text{for all } k \in \mathcal{K}_S \text{ with } z_k^* > 0. \tag{B.27}$$

From (B.21) and (B.25), we have for large enough $r$,

$$\sum_{i\in\mathcal{I}} y_i^r R_{ij}^r = \sum_{k\in\mathcal{K}_S} A_{kj} z_k^r \qquad \text{for all } j \in \mathcal{J}_S \text{ with } a_j^* > 0.$$

Therefore,

$$\sum_{j\in\mathcal{J}_S} a_j^* \sum_{i\in\mathcal{I}} y_i^r R_{ij}^r = \sum_{j\in\mathcal{J}_S} a_j^* \sum_{k\in\mathcal{K}_S} A_{kj} z_k^r \qquad \text{for large enough } r. \tag{B.28}$$

Similarly, it follows from (B.22) and (B.26) that, for large enough $r$,

$$\sum_{j\in\mathcal{J}_I} a_j^* \sum_{i\in\mathcal{I}} y_i^r R_{ij}^r = -\sum_{j\in\mathcal{J}_I} a_j^* \sum_{k\in\mathcal{K}_I} A_{kj} z_k^r. \tag{B.29}$$

From (B.23) and (B.27), it follows that, for large enough $r$,

$$\sum_{j\in\mathcal{J}_s} a_j^* A_{kj} = 1 \qquad \text{for all } k \in \mathcal{K}_S \text{ with } z_k^r > 0.$$

Therefore,

$$\sum_{k\in\mathcal{K}_S} z_k^r \sum_{j\in\mathcal{J}_S} a_j^* A_{kj} = \sum_{k\in\mathcal{K}_S} z_k^r \qquad \text{for large enough } r. \tag{B.30}$$

It follows from (B.28) and (B.29) that

$$\begin{aligned}\sum_{j\in\mathcal{J}} a_j^* \sum_{i\in\mathcal{I}} y_i^r R_{ij}^r &= \sum_{j\in\mathcal{J}} a_j^* \sum_{k\in\mathcal{K}_S} A_{kj} z_k^r - \sum_{j\in\mathcal{J}} a_j^* \sum_{k\in\mathcal{K}_I} A_{kj} z_k^r \\ &= 1 - \rho^r \qquad \text{for large enough } r.\end{aligned}$$

The second equality follows from (B.30), $\sum_{k\in\mathcal{K}_S} z_k^r = 1$, $\sum_{j\in\mathcal{J}_I} a_j^* A_{kj} = 1$ for each $k \in \mathcal{K}_I$, and $\sum_{k\in\mathcal{K}_I} z_k^r = \rho^r$.

Then Lemma 10 follows from the fact that $\max_{a\in\mathcal{E}} y^r \cdot R^r a = 1 - \rho^r$. □

PROOF OF THEOREM 1. We define the scaled process $\mathbb{X}^r$ via

$$\mathbb{X}^r(t) = r^{-2}\mathbb{X}(r^2 t) \qquad \text{for each } t \geq 0.$$

Fix a sample path that satisfies the strong law of large numbers for $u_j$ and $\phi_j^i$. Let $(\bar{\bar{Z}}, \bar{\bar{T}})$ be a fluid limit of $(\bar{\bar{Z}}^r, \bar{\bar{T}}^r)$ along the sample path. Following the arguments in Section A.2 in Dai and Lin (2005), such a limit exists and satisfies the fluid model equations (6.1)–(6.5) presented in Section 6. Under the maximum pressure policy with parameter $\alpha$, each fluid limit $(\bar{\bar{Z}}, \bar{\bar{T}})$ also satisfies the fluid model equation (6.6). The justification of fluid model equation (6.6) is similar to Lemma 4 in Dai and Lin (2005), with the scaling $r$ replaced by $r^2$. Therefore, $(\bar{\bar{Z}}^r, \bar{\bar{T}}^r)$ is a fluid model solution under the maximum pressure policy. Similar to the proof of Theorem 4 in Dai and Lin (2005), using a different Lyapunov function $f(t) = (\alpha \times Z(t)) \cdot Z(t)$, one can easily prove that the fluid model under the maximum pressure policy with a general parameter set $(\alpha, \beta)$ is weakly stable; namely, any fluid model solution $(\bar{\bar{Z}}^r, \bar{\bar{T}}^r)$ under the maximum pressure policy satisfies $Z(t) = 0$ for each $t \geq 0$ given $Z(0) = 0$. As a consequence, we have for any $t > 0$, $\bar{\bar{T}}(t)/t$ satisfies (3.2)–(3.5) with $\rho = 1$. Because $x^*$ is the unique optimal solution to the static planning problem (3.1)–(3.5) with objective value equal to 1, $\bar{\bar{T}}(t) = xt$ for each $t \geq 0$. Since this is true for any fluid limit, we have $\bar{\bar{T}}^r(t) \to x^* t$ for each $t$ with probability 1, which implies asymptotic efficiency. □

**Acknowledgments.** The authors thank two anonymous referees for significantly improving the presentation of the paper.

H. Milton Stewart School of Industrial and Systems Engineering
Georgia Institute of Technology
Atlanta, Georgia 30332
USA
E-mail: dai@gatech.edu

Kellogg School of Management
Northwestern University
Evanston, Illinois 60208
USA
E-mail: wuqin-lin@kellogg.northwestern.edu